# Weighted cohomology of arithmetic groups

By Arvind Nair*

## Introduction

Let $\mathbf{G}$ be a semisimple algebraic group defined over the rational numbers, $K$ a maximal compact subgroup of $G = \mathbf{G}(\mathbb{R})$, and $\Gamma \subset \mathbf{G}(\mathbb{Q})$ a neat arithmetic subgroup. Let $X = \Gamma \backslash G/K$ be the locally symmetric space associated to $\Gamma$, and $\mathbb{E}$ the local system on $X$ constructed out of a finite-dimensional, irreducible, algebraic representation $E$ of $\mathbf{G}$. Fix a maximally $\mathbb{Q}$-split torus $\mathbf{S}$ in $\mathbf{G}$; $\mathbf{S}$ is assumed to be nontrivial, so that $X$ is noncompact. Let $A = \mathbf{S}(\mathbb{R})^0$ and, by a slight abuse of notation, let $X^*(A)$ denote the group of rational characters of $\mathbf{S}$. Choose a minimal rational parabolic subgroup $\mathbf{P}_0 \supset \mathbf{S}$; the choice determines a notion of positivity in $X^*(A)$ and a set of positive roots among the roots of $\mathbf{S}$ in $\mathbf{G}$. Let $\rho_0$ be the half-sum of the positive roots.

*Weighted cohomology* is an invariant of $\Gamma$ introduced by M. Goresky, G. Harder, and R. MacPherson [12] in the study ([15], [16], [12], [14]) of the trace of Hecke operators in the cohomology of $\Gamma$. For each $p \in X^*(A) \otimes \mathbb{Q}$ there are *weight profiles* $\overline{p}$ and $\underline{p}$ (*upper* and *lower* $p$; see 1.5) and groups $W^{\overline{p}}H^i(\Gamma, E)$ and $W^{\underline{p}}H^i(\Gamma, E)$. (The notion of profile here differs from (but is equivalent to) that of [12]; see Remark 1.2 below.) If $p+q = -2\rho_0$ then the profiles $\overline{p}$ and $\underline{q}$ are dual and the corresponding groups are in Poincaré duality. For very positive $p$ one gets the compactly supported cohomology $H^i_c(X, \mathbb{E})$ and dually, for very negative $p$, the full cohomology $H^i(X, \mathbb{E})$. The definition of the weighted cohomology groups is somewhat involved, so I will attempt to give some idea of the motivation behind it; detailed definitions are recalled in Section 1. Associated to $\overline{p}$ and $\underline{p}$ are complexes of sheaves $\mathbf{W}^{\overline{p}}\mathbf{C}^\bullet(\mathbb{E})$ and $\mathbf{W}^{\underline{p}}\mathbf{C}^\bullet(\mathbb{E})$ on the reductive Borel-Serre compactification $\overline{X}$ of $X$. Let $\mathbf{P}$ be a rational parabolic subgroup, $\mathbf{N}_\mathbf{P}$ its unipotent radical, $\mathbf{S}_\mathbf{P}$ the split centre of $\mathbf{P}/\mathbf{N}_\mathbf{P}$ and $X_P$ the $\mathbf{P}$-boundary component in $\overline{X}$. The idea behind the construction in [12] is as follows: The stalk cohomology of the direct image $Ri_*\mathbb{E}$ (by $i: X \hookrightarrow \overline{X}$) at $x \in X_P$ is given by the Lie algebra cohomology $H^*_{\mathfrak{N}_P}(E)$ (here $\mathfrak{N}_P = \text{Lie}(\mathbf{N_P}(\mathbb{R})))$. This has a natural action of $\mathbf{S_P}(\mathbb{R})$, which is realized geometrically for elements in $\mathbf{S_P}(\mathbb{Q})$ by Looijenga's Hecke correspondences ([20, 3.7], [12, 15]). The profile $\overline{p}$ defines a canonical submodule $H^*_{\mathfrak{N}_P}(E)_+ \subset H^*_{\mathfrak{N}_P}(E)$ with "high" $\mathbf{S_P}(\mathbb{R})$-weights (1.5). The idea is then to truncate (see [12, 1.8]) $Ri_*\mathbb{E}$ along the boundary to get an object $\mathbf{W}^{\overline{p}}\mathbf{C}^\bullet(\mathbb{E})$ with stalk cohomology $H^*_{\mathfrak{N}_P}(E)_+$ at $x$. The difficulty is in finding a version of $Ri_*\mathbb{E}$ with a semisimple action of $\mathbf{S_P}(\mathbb{Q})$ representing the geometric action. (Two constructions of $\mathbf{W}^{\overline{p}}\mathbf{C}^\bullet(\mathbb{E})$ are given in [12]; the one used in this paper is recalled in 1.4, 1.5.) Then $W^{\overline{p}}H^i(\Gamma, E)$ is the $i^{\text{th}}$ hypercohomology group of $\mathbf{W}^{\overline{p}}\mathbf{C}^\bullet(\mathbb{E})$. The profile $\underline{p}$ defines a slightly larger submodule $H^*_{\mathfrak{N}_P}(E)_-$ (see 1.5) and then the same construction defines $\mathbf{W}^{\underline{p}}\mathbf{C}^\bullet(\mathbb{E})$ and $W^{\underline{p}}H^i(\Gamma, E)$.

*Supported by a Sloan Foundation doctoral dissertation fellowship at the University of Michigan and by NSF grant DMS 9304580 at the Institute for Advanced Study.

The space $X$ is a Riemannian manifold, a complete metric coming from an invariant metric on $G/K$. There is a natural choice of inner product (1.2) on $E$ metrizing $\mathbb{E}$. For each $\lambda \in \mathfrak{A}^* = X^*(A) \otimes \mathbb{R}$, J. Franke [11] has defined *weighted $L^2$ cohomology groups* via complexes of forms on $X$ satisfying certain square-integrability conditions given by $\lambda$. For very positive $\lambda$ these give the full cohomology $H^i(X, \mathbb{E})$ and for very negative $\lambda$ they give $H^i_c(X, \mathbb{E})$. The definition of these groups involves a little reduction theory. Recall that a Siegel set in $G$ (with respect to $\mathbf{P}_0 = \mathbf{NMS}$) is a set of the form $\mathfrak{S} = \Omega.A(t).K$ where $\Omega$ is a relatively compact set in $\mathbf{N}(\mathbb{R})\mathbf{M}(\mathbb{R})$ and $A(t) \subset A$ is the set of $a$ such that $a^\alpha > t$ for every simple root $\alpha$ of $\mathbf{S}$ in $\mathbf{G}$. Reduction theory [2] says that, for large enough $\Omega$ and $t$, there is a finite set $C \subset \mathbf{G}(\mathbb{Q})$ such that $C.\mathfrak{S}$ is a coarse fundamental domain for $\Gamma$ in $G$. Now for $\lambda \in \mathfrak{A}^*$, there is an *admissible weight function* $w_\lambda$ on $G$ which descends to $X$ (see 1.7 for the precise definition). The essential property of $w_\lambda$ is that both $a^\lambda w_\lambda(cnmak)$ and $a^{-\lambda}w_\lambda(cnmak)^{-1}$ are bounded functions of $cnmak \in C.\mathfrak{S}$. Informally, one might say that $w_\lambda$ grows like $a^{-\lambda}$ as one goes to the boundary of $X$. Consider the space $L^i_{\lambda-\log}(X, \mathbb{E})$ of measurable $\mathbb{E}$-valued $i$-forms $\sigma$ on $X$ such that

$$(*) \qquad w_\lambda \log(w_\epsilon)^j \, \sigma \text{ and } w_\lambda \log(w_\epsilon)^j \, d\sigma \text{ are square-integrable on } X$$

for all $j \geq 0$ (here $d$ is the distributional exterior derivative and $\epsilon$ is some fixed dominant weight on $A$). The admissibility of $w_\lambda$ (see 1.7) ensures that $(L^\bullet_{\lambda-\log}(X, \mathbb{E}), d)$ is a complex; define $H^i_{\lambda-\log}(X, \mathbb{E})$ to be its $i^{\text{th}}$ cohomology group. Requiring $(*)$ to hold for some $j \leq 0$ (instead of all $j \geq 0$) gives a complex $L^\bullet_{\lambda+\log}(X, \mathbb{E}) \supset L^\bullet_{\lambda-\log}(X, \mathbb{E})$ and its cohomology is denoted $H^i_{\lambda+\log}(X, \mathbb{E})$. For $\lambda = 0$ the groups will be denoted $H^i_{(2)-\log}(X, \mathbb{E})$ and $H^i_{(2)+\log}(X, \mathbb{E})$.

The main theorem is that these two cohomology theories are the same:

THEOREM A. *Let $\lambda = -p - \rho_0$. There are natural isomorphisms*
$$W^{\overline{p}}H^i(\Gamma, E) \simeq H^i_{\lambda-\log}(X, \mathbb{E}),$$
$$W^{\underline{p}}H^i(\Gamma, E) \simeq H^i_{\lambda+\log}(X, \mathbb{E}).$$

The two assertions are equivalent by duality (see Remarks 1.3 and 1.9).

Consider the space $S_{\lambda-\log}(\Gamma\backslash G)$ (resp. $S_{\lambda+\log}(\Gamma\backslash G)$) of smooth complex-valued functions $f$ on $\Gamma\backslash G$ such that, for all $j \geq 0$ (resp. some $j \leq 0$),

$$w_\lambda \log(w_\epsilon)^j R_D f \text{ is square-integrable}$$

for every operator $D$ in the universal enveloping algebra of $\mathfrak{G} = \text{Lie}(G)$. (By a Sobolev lemma (see [11, Prop. 2.3.2]), $S_{\lambda-\log}(\Gamma\backslash G)$ is also the space of functions such that, for all $j \geq 0$, $w_{\lambda+\rho_0}\log(w_\epsilon)^j R_D f$ is bounded for all $D$ and similarly for $S_{\lambda+\log}(\Gamma\backslash G)$.) They are $(\mathfrak{G}, K)$-modules. A theorem of Franke's (see Proposition 1.4 below) computes $H^i_{\lambda\pm\log}(X, \mathbb{E})$ in terms of $S_{\lambda\pm\log}(\Gamma\backslash G)$, and so Theorem A is equivalent to:

COROLLARY A. *Let $\lambda = -p - \rho_0$. There are natural isomorphisms*
$$W^{\overline{p}}H^i(\Gamma, E) \simeq H^i_{(\mathfrak{G},K)}(S_{\lambda-\log}(\Gamma\backslash G) \otimes E),$$
$$W^{\underline{p}}H^i(\Gamma, E) \simeq H^i_{(\mathfrak{G},K)}(S_{\lambda+\log}(\Gamma\backslash G) \otimes E).$$

There are two special profiles, *upper middle* $\mu = \overline{-\rho_0}$ and *lower middle* $\nu = \underline{-\rho_0}$. They are dual to each other. Theorem A says that $W^\mu H^i(\Gamma, E) \simeq H^i_{(2)-\log}(X, \mathbb{E})$ and $W^\nu H^i(\Gamma, E) \simeq H^i_{(2)+\log}(X, \mathbb{E})$. The relation of $H^i_{(2)\pm\log}(X, \mathbb{E})$ to the usual $L^2$ cohomology is given by:

THEOREM B. *If the ranks of $G$ and $K$ are equal then*
$$H^i_{(2)-\log}(X, \mathbb{E}) \simeq H^i_{(2)}(X, \mathbb{E}) \simeq H^i_{(2)+\log}(X, \mathbb{E}).$$

In fact this holds under a slightly more general condition from [5] (see Theorem 4.1). Theorems A and B imply:

COROLLARY B. *If the ranks of $G$ and $K$ are equal then*
$$W^\mu H^i(\Gamma, E) \simeq H^i_{(2)}(X, \mathbb{E}) \simeq W^\nu H^i(\Gamma, E).$$

This corollary is related to Zucker's conjecture (a theorem thanks to Looijenga [20] and Saper-Stern [24]). Suppose that $X$ is Hermitian and $\hat{X}$ is its Baily-Borel Satake compactification; $\hat{X}$ is a complex projective variety. A theorem of Goresky-Harder-MacPherson ([12, Th. 23.2]) says that $W^\mu H^i(\Gamma, E)$ and $W^\nu H^i(\Gamma, E)$ are isomorphic to the intersection cohomology $IH^i(\hat{X}, \mathbb{E})$.

The corollary then implies Zucker's conjecture. (Note that the proof of Theorem 23.2 of [12] uses the main technical step of [20], so this is not a new proof of the conjecture.)

The proof of Theorem A has three steps. The first step is the equivalence with Corollary A (due to Franke) and is dealt with in Section 1. In this section, I also recall the construction of $\mathbf{W}^{\overline{p}}\mathbf{C}^\bullet(\mathbb{E})$ and $\mathbf{W}^{\underline{p}}\mathbf{C}^\bullet(\mathbb{E})$ and introduce complexes of fine sheaves $\mathbf{S}^\bullet_{\lambda\pm\log}(\mathbb{E})$ on $\overline{X}$. The complex of global sections of $\mathbf{S}^\bullet_{\lambda\pm\log}(\mathbb{E})$ is the standard Lie algebra cohomology complex computing $H^*_{(\mathfrak{G},K)}(S_{\lambda\pm\log}(\Gamma\backslash G)\otimes E)$. So to prove Corollary A it suffices to prove that $\mathbf{W}^{\overline{p}}\mathbf{C}^\bullet(\mathbb{E})$ and $\mathbf{S}^\bullet_{\lambda-\log}(\mathbb{E})$ are quasi-isomorphic (the assertion about $\underline{p}$ follows by duality). The following remarks motivate this, without actually playing a role in the proof: The complex of sheaves $\mathbf{S}^\bullet_\infty(\mathbb{E}) = \varinjlim_\lambda \mathbf{S}^\bullet_{\lambda-\log}(\mathbb{E})$ is a version of $Ri_*\mathbb{E}$, and the subsheaf defined by the condition $(*)$ (which can be made local) is $\mathbf{S}^\bullet_{\lambda-\log}(\mathbb{E})$. The logarithmic terms in $(*)$ make $\mathbf{S}^\bullet_{\lambda-\log}(\mathbb{E})$ cohomologically constructible for the natural stratification of $\overline{X}$. The stalk cohomology of $\mathbf{S}^\bullet_{\lambda-\log}(\mathbb{E})$ is a submodule of $H^*_{\mathfrak{N}_P}(E)$ and for $\lambda$ chosen correctly (i.e. $\lambda = -p - \rho_0$) it is exactly the submodule $H^*_{\mathfrak{N}_P}(E)_+$ defined by $\overline{p}$. (For $\lambda = 0$ this is essentially in [29].) So, like $\mathbf{W}^{\overline{p}}\mathbf{C}^\bullet(\mathbb{E})$, $\mathbf{S}^\bullet_{\lambda-\log}(\mathbb{E})$ is constructed from a version of $Ri_*\mathbb{E}$ by a truncation procedure and has the same local cohomological properties, suggesting that these two should be quasi-isomorphic. In fact, there is an explicit and natural quasi-isomorphism, as shown in the next two steps.

The second step (§2) in the proof of Theorem A is to show that there is an inclusion $\mathbf{W}^{\overline{p}}\mathbf{C}^\bullet(\mathbb{E}) \hookrightarrow \mathbf{S}^\bullet_{\lambda-\log}(\mathbb{E})$ for $\lambda = -p - \rho_0$. The final step (§3) is to show that this inclusion is a quasi-isomorphism via a local calculation on $\overline{X}$. There are two key points here: Proposition 3.2 (reduction to the space of $N_P$-invariants of $\mathbf{S}^\bullet_{\lambda-\log}(\mathbb{E})$ over neighbourhoods of boundary points, (a version of) an idea that goes back to van Est [27]) and Lemma 3.5 (where the logarithmic terms in $(*)$ become essential). In Section 0 below, for the reader's convenience, the proof of Theorem A is sketched in the simplest possible case, namely when $\mathbf{G} = \mathrm{SL}(2)$ and $E$ is trivial.

The proof of Theorem B, which is global and independent of the previous sections, is in Section 4 and uses results from [11] and [5] and standard Lie algebra cohomology arguments. The final section (§5) contains applications and several comments on relations to other work and on generalizations.

*Acknowledgements.* These results were part of my doctoral dissertation [21] at the University of Michigan. It is a great pleasure to thank Professor Gopal Prasad for his constant support and encouragement and for many helpful comments. I thank Professors A. Borel, W. Casselman, M. Goresky, and R. MacPherson for their comments on earlier versions of this work and the Institute for Advanced Study for its hospitality and support while this paper was substantially revised. I am indebted to A. Borel for numerous detailed comments and for help with the revision, to W. Casselman for suggesting the more group-theoretic approach (following [9]) used in this revised version and to M. Goresky for several helpful conversations and for improvements in the exposition. Finally, I am grateful to the referees for helpful comments and corrections, in particular for a simplification in the proof of Theorem 4.1.

## 0. Example: $\mathbf{G} = \mathrm{SL}(2), E = \mathbb{C}$

The proof here is complete modulo one or two technical details for which easy proofs can be supplied in this case and which are, in any case, treated in general later. Note that 0.n corresponds to §n in the body of the paper.

0.1. Let $\mathbf{G} = \mathrm{SL}(2)$, $K = \mathrm{SO}(2)$, and identify $D = \mathbf{G}(\mathbb{R})/K$ with the upper half-plane with the coordinates $x, y$. Assume that $E = \mathbb{C}$ is the trivial representation. The reductive Borel-Serre compactification adds finitely many cusps to the modular curve $X = \Gamma\backslash D$. I shall assume, for simplicity, that $\Gamma \subset \mathrm{SL}(2,\mathbb{Z})$ has only one cusp, *viz.* $i\infty$. In the Borel-Serre compactification the cusp is blown up to a circle at $y = \infty$. The stabilizer of $i\infty$ is the Borel subgroup $B$ of upper triangular matrices. Let $N$ be the upper triangular unipotent group and $\Gamma_N = \Gamma \cap N$. Let $A = \{\mathrm{diag}(a, a^{-1}) \,|\, a > 0\}$ and let $X^*(A)$ be its group of algebraic characters. Since $X^*(A)$ has a canonical generator $\beta$ defined by $\beta(\mathrm{diag}(a, a^{-1})) = a^2$, we shall identify it with $\mathbb{Z}\beta$. Let $\mathrm{pr}_\Gamma : D \to X$ and $\mathrm{pr}_K : \Gamma\backslash G \to X$ be the obvious projections. A fundamental system of neighbourhoods of $i\infty$ in $\overline{X}$ is given by

$$V_0(T) := \mathrm{pr}_\Gamma(\{(x,y) \,|\, y > T\}) \cup \{i\infty\}.$$

There is a retraction $r_{i\infty} : D \to \{y = \infty\}$ given by $r_{i\infty}(x,y) = (x, y = \infty)$.

A differential form $\sigma$ on $X$ is *special* if there is a neighbourhood $V_0(T)$ on which $\sigma$ is the pullback under $r_{i\infty}$ of a translation-invariant (i.e., $N$-invariant) form $\sigma_{i\infty}$ on the line $y = \infty$. The special forms give a presheaf on $\overline{X}$ in which the sections over $U \subset \overline{X}$ are the restrictions of special forms on $X$ to $U \cap X$. Translation-invariant forms on $\{y = \infty\}$ are identified with $\Lambda^\bullet \mathfrak{N}^*$, which has a natural action of $A$. Let $p \in \mathbb{Q}$. Define a presheaf $P\mathbf{W}^{\overline{p}}\mathbf{C}^i(\mathbb{C})$ by letting its sections over $U$ be the forms that, for some $T$, are special on $U_{0,T} = U \cap V_0(T)$ and are lifted from an element of $\Lambda^\bullet \mathfrak{N}^*$ with $A$-weight strictly greater than $p\beta$. Let $\mathbf{W}^{\overline{p}}\mathbf{C}^i(\mathbb{C})$ be the sheaf associated to this presheaf. If $U \subset X$ then $\mathbf{W}^{\overline{p}}\mathbf{C}^\bullet(U, \mathbb{C}) = \Omega^\bullet(U, \mathbb{C})$. The complex of sheaves $\mathbf{W}^{\overline{p}}\mathbf{C}^\bullet(\mathbb{C})$ is fine and $W^{\overline{p}}H^i(\Gamma, \mathbb{C})$ is its $i^{\mathrm{th}}$ cohomology group.

As $A$-modules, $\Lambda^0\mathfrak{N}^* = \mathbb{C}$ and $\Lambda^1\mathfrak{N}^* = \mathfrak{N}^* = \mathbb{C}_{-\beta}$. Varying $p$, one gets three complexes of sheaves up to isomorphism:

$$\mathbf{W}^{\infty}\mathbf{C}^{\bullet}(\mathbb{C}) = i_!\Omega^{\bullet}(\mathbb{C}), \quad \mathbf{W}^{\overline{-\frac{1}{2}}}\mathbf{C}^{\bullet}(\mathbb{C}), \quad \mathbf{W}^{-\infty}\mathbf{C}^{\bullet}(\mathbb{C}) \simeq Ri_*\mathbb{C}.$$

(In this case $\hat{X} = \overline{X}$ and and the middle weighted cohomology complex $\mathbf{W}^{\overline{-\frac{1}{2}}}\mathbf{C}^{\bullet}(\mathbb{C})$ is quasi-isomorphic to the intersection complex on $\hat{X}$. This is the trivial case of Theorem 23.2 of [12].)

Let $\lambda\beta \in \mathfrak{A}^* = \mathbb{R}\beta$. For $D \in \mathcal{U}(\mathfrak{G})$ let $R_D$ be the associated invariant differential operator on functions on $\Gamma\backslash G$. For $W \subset \Gamma\backslash G$ and stable under $K$, let $S_{\lambda-\log}(W)$ be the space of smooth $K$-finite functions $f : W \to \mathbb{C}$ such that, for all $m \geq 0$,

$$\log(y)^m e^{-\lambda\log(y)} R_D f \in L^2(W)$$

for all $D \in \mathcal{U}(\mathfrak{G})$. It is a $(\mathfrak{G}, K)$-module. Define a presheaf $P\mathbf{S}^i_{\lambda-\log}(\mathbb{C})$ on $\overline{X}$ by

$$P\mathbf{S}^i_{\lambda-\log}(U, \mathbb{C}) = \mathrm{Hom}_K\big(\Lambda^i(\mathfrak{G}/\mathfrak{K}), S_{\lambda-\log}(\mathrm{pr}_K^{-1}(U \cap X))\big).$$

The associated sheaf is denoted $\mathbf{S}^i_{\lambda-\log}(\mathbb{C})$ and it is fine. For $U \subset \overline{X}$, $P\mathbf{S}^i_{\lambda-\log}(U,\mathbb{C})$ is identified with a subspace of the smooth $i$-forms on $U \cap X$ (see 2.1). The cohomology groups of the complex of sheaves $\mathbf{S}^{\bullet}_{\lambda-\log}(\mathbb{C})$ are $H^*_{\lambda-\log}(X,\mathbb{C})$. (Results of Franke's recalled in 1.8 imply that this agrees with the definition given earlier in the introduction.)

0.2. I want to show that there is an inclusion $\mathbf{W}^{\overline{p}}\mathbf{C}^i(\mathbb{C}) \to \mathbf{S}^i_{\lambda-\log}(\mathbb{C})$ when $\lambda \geq -p - 1/2$ (the second step referred to in the introduction). In fact this holds for the presheaves used to define these sheaves.

Suppose $U \subset \overline{X}$ is open and $\sigma \in P\mathbf{W}^{\overline{p}}\mathbf{C}^i(U,\mathbb{C})$. By definition, $\sigma$ is the restriction to $U$ of a special form on $X$. To show that $\sigma \in P\mathbf{S}^i_{\lambda-\log}(U,\mathbb{C})$ it will suffice to show that $\sigma \in P\mathbf{S}^i_{\lambda-\log}(V_0(t_0),\mathbb{C})$ for some $t_0$, to which end we need a description of this last space.

Define $S_{\lambda-\log}(\Gamma_N\backslash NA(t_0))$ to be the space of smooth functions $f$ on $\Gamma_N\backslash NA(t_0)$ such that, for any $m \geq 0$, $\log(y)^m e^{-\lambda\log(y)} R_D f$ is, for any $D \in \mathcal{U}(\mathfrak{B})$, square-integrable with respect to the measure $y^{-2}\,dx\,dy$ coming from the left-invariant measure on $G$ (here $NA(t_0)$ is identified with $\{y > t_0^{1/2}\} \subset D$, so we use the coordinates $x, y$). There is an isomorphism

$$\mathrm{Hom}_K(\Lambda^{\bullet}(\mathfrak{G}/\mathfrak{K}), S_{\lambda-\log}(\Gamma_N\backslash NA(t_0)K)) \cong \mathrm{Hom}(\Lambda^{\bullet}\mathfrak{B}, S_{\lambda-\log}(\Gamma_N\backslash NA(t_0)))$$

because $S_{\lambda-\log}(\Gamma_N\backslash NA(t_0)K) = \mathrm{Ind}_{(\mathfrak{B},e)}^{(\mathfrak{G},K)}(S_{\lambda-\log}(\Gamma_N\backslash NA(t_0)))$ (for the appropriate induction functor; see Lemma 2.1). For $\delta\beta \in \mathbb{R}\beta$, let $S_{\delta-\log}(A(t))$ be the space of smooth functions $f$ on $A(t)$ such that $\log(y)^m e^{-\delta\log(y)} R_D f \in L^2(A(t), dy/y)$ for any $m \geq 0$, for any $D \in \mathcal{U}(\mathfrak{A})$. Restricting attention to

$N$-invariants we have
$$S_{\lambda-\log}(\Gamma_N\backslash NA(t_0))^{N-\mathrm{inv}} \cong S_{\lambda+1/2-\log}(A(t)).$$

(To see this, note that the algebra $\{R_D|D \in \mathcal{U}(\mathfrak{B})\}$ is generated by $y\,\partial/\partial y$ and $y\,\partial/\partial x$, so if $f$ is $N$-invariant then only the condition with respect to powers of $y\,\partial/\partial y$ matters.) Via this isomorphism, the complex of $N$-invariant sections of $P\mathbf{S}^\bullet_{\lambda-\log}(\mathbb{C})$ over $V_0(t)$ is identified with
$$\mathrm{Hom}(\Lambda^\bullet\mathfrak{A} \otimes \Lambda^\bullet\mathfrak{N}, S_{\lambda+1/2-\log}(A(t))).$$

This is identified with the double complex
$$\mathrm{Hom}\bigl(\Lambda^\bullet\mathfrak{A}, S_{\lambda+1/2-\log}(A(t)) \otimes \Lambda^\bullet\mathfrak{N}^*\bigr)$$

(see Lemma 2.3). Now $\Lambda^\bullet\mathfrak{N}^* = \oplus_\alpha (\Lambda^\bullet\mathfrak{N}^*)_\alpha$ as a complex (the differential commutes with the $A$-action) and so
$$P\mathbf{S}^\bullet_{\lambda-\log}(V_0(t),\mathbb{C})^{N-\mathrm{inv}} \cong \bigoplus_\alpha \mathrm{Hom}\bigl(\Lambda^\bullet\mathfrak{A}, S_{\lambda+1/2-\log}(A(t)) \otimes (\Lambda^\bullet\mathfrak{N}^*)_\alpha\bigr).$$

Let $V_\alpha^\bullet$ be the complex $(\Lambda^\bullet\mathfrak{N}^*)_\alpha$ considered as a trivial $A$-module and let $\mathbb{C}_\alpha$ be the one-dimensional $A$-module of character $\alpha$. Then, as $A$-modules, $(\Lambda^\bullet\mathfrak{N}^*)_\alpha \cong \mathbb{C}_\alpha \otimes V_\alpha^\bullet$ canonically. Furthermore, there is an isomorphism of $\mathfrak{A}$-modules
$$S_{\lambda+1/2-\log}(A(t)) \otimes \mathbb{C}_\alpha \cong S_{\lambda+1/2+\alpha-\log}(A(t))$$

(this is easily checked; see 2.4). Combining these, we have
$$P\mathbf{S}^\bullet_{\lambda-\log}(V_0(t),\mathbb{C})^{N-\mathrm{inv}} \cong \bigoplus_\alpha \bigl(\Lambda^\bullet\mathfrak{A}^* \otimes S_{\lambda+1/2+\alpha-\log}(A(t))\bigr) \otimes V_\alpha^\bullet$$

as a complex.

Let $A_\alpha$ be the space of constant functions on $A(t)$ if $\alpha > p$ and zero otherwise. Then $A_\alpha \subset S_{\lambda+1/2+\alpha-\log}(A(t))$ since $\lambda \geq -p - 1/2$ and so
$$W_p^\bullet(t) = \bigoplus_\alpha \Lambda^0\mathfrak{A}^* \otimes A_\alpha \otimes V_\alpha^\bullet$$

with differential $\mathrm{id} \otimes (d : V_\alpha^i \to V_\alpha^{i+1})$ is a subcomplex of $P\mathbf{S}^\bullet_{\lambda-\log}(V_0(t),\mathbb{C})^{N-\mathrm{inv}}$. Now if $\sigma$ is lifted from the boundary on $V_0(t_0)$ then $\phi_\sigma$ belongs to $W_p^i(t)$ and hence to $P\mathbf{S}^\bullet_{\lambda-\log}(V_0(t_0),\mathbb{C})^{N-\mathrm{inv}}$. Therefore there is an injection $P\mathbf{W}^{\overline{p}}\mathbf{C}^i(\mathbb{C}) \to P\mathbf{S}^i_{\lambda-\log}(\mathbb{C})$ when $\lambda \geq -p - 1/2$.

0.3. I want to show that the inclusion $\mathbf{W}^{\overline{p}}\mathbf{C}^\bullet(\mathbb{C}) \to \mathbf{S}^\bullet_{\lambda-\log}(\mathbb{C})$ gives an isomorphism on stalk cohomology at the cusp when $\lambda = -p - 1/2$. This is the third step in the proof of Theorem A.

For $W \subset \Gamma\backslash G$ such that $\mathrm{pr}_\Gamma^{-1}(W)$ is left $N$-stable, let $S_{\mathrm{csp}}(W)$ be the space of functions $f \in S_{\lambda-\log}(W)$ such that the constant term of $\sigma$ along $B$ vanishes on the intersection of $W$ and $\Gamma \backslash NA(t)K = \mathrm{pr}_K^{-1}(V_0(t) \cap X)$ for some $t$. When $W$ is right $K$-stable this is a $(\mathfrak{G}, K)$-module. (Note that the functions in $S_{\mathrm{csp}}(\Gamma \backslash G)$ are not necessarily cuspidal.) The constant term along $B$ gives a surjective map of $(\mathfrak{G}, K)$-modules

$$S_{\lambda-\log}(\Gamma_N \backslash NA(t_0)K) \;\to\; \varinjlim_{t > t_0} S_{\lambda-\log}(\Gamma_N \backslash NA(t)K)^{N-\mathrm{inv}}$$

with kernel $S_{\mathrm{csp}}(\Gamma_N \backslash NA(t_0)K)$ and so a surjection of complexes

$$C^\bullet_{(\mathfrak{G},K)}(\varinjlim S_{\lambda-\log}(\Gamma_N \backslash NA(t)K)) \;\xrightarrow{\mathrm{CT}}\; C^\bullet_{(\mathfrak{G},K)}(\varinjlim S_{\lambda-\log}(\Gamma_N \backslash NA(t)K)^{N-\mathrm{inv}})$$

(direct limits are as $t$ grows). Note that, by definition, the domain of CT is the stalk of $\mathbf{S}^\bullet_{\lambda-\log}(\mathbb{C})$ at $i\infty$ and its range is the $N$-invariants in the stalk. The kernel of CT is the complex

$$C^\bullet_{(\mathfrak{G},K)}(\varinjlim S_{\mathrm{csp}}(\Gamma_N \backslash NA(t)K)).$$

This is an induced module and so has cohomology

$$H^*_{(\mathfrak{G},K)}(\varinjlim S_{\mathrm{csp}}(\Gamma_N \backslash NA(t)K)) \cong H^*_{\mathfrak{A}+\mathfrak{N}}(\varinjlim S_{\mathrm{csp}}(\Gamma_N \backslash NA(t)))$$

(see proof of Proposition 3.2). The Hochschild-Serre spectral sequence converging to this has $E_2$-term

$$H^*_{\mathfrak{A}}(H^*_{\mathfrak{N}}(\varinjlim S_{\mathrm{csp}}(\Gamma_N \backslash NA(t)))).$$

In fact the $\mathfrak{N}$-cohomology groups vanish by a well-known argument using Fourier series on the circle $\Gamma_N \backslash N$ (see [19, 2.4]). Therefore CT is a quasi-isomorphism.

The complex $W^\bullet_p(t) = \bigoplus_\alpha \Lambda^0 \mathfrak{A}^* \otimes A_\alpha \otimes V^\bullet_\alpha$ of 0.2 is identified with the space of forms on $V_0(t)$ that are lifted from a form with weight $> p$ on the cusp. Given an element $\sigma \in P\mathbf{W}^{\overline{p}}\mathbf{C}^i(V_0(t_0), \mathbb{C})$ its constant term belongs to $W^i_p(t)$ for some $t > t_0$. This defines a map

$$P\mathbf{W}^{\overline{p}}\mathbf{C}^i(V_0(t_0), \mathbb{C}) \;\to\; \varinjlim W^i_p(t)$$

with kernel the space of $i$-forms vanishing near the cusp. The map is surjective, and so the constant term gives an identification of the stalk of $P\mathbf{W}^{\overline{p}}\mathbf{C}^\bullet(\mathbb{C})$ at $i\infty$ with $\varinjlim W^i_p(t)$. Since CT is a quasi-isomorphism, we see that $\mathbf{W}^{\overline{p}}\mathbf{C}^\bullet(\mathbb{C}) \to \mathbf{S}^\bullet_{\lambda-\log}(\mathbb{C})$ is a quasi-isomorphism at $i\infty$ if and only if the inclusion

$$\varinjlim W^\bullet_p(t) \;\to\; \varinjlim P\mathbf{S}^\bullet_{\lambda-\log}(V_0(t), \mathbb{C})^{N-\mathrm{inv}}$$

is a quasi-isomorphism. As in 0.2, the right-hand side is

$$\bigoplus_\alpha (\Lambda^\bullet \mathfrak{A}^* \otimes S_{\lambda+1/2+\alpha-\log}(A(t))) \otimes V^\bullet_\alpha.$$

It suffices to show that the inclusion $A_\alpha \hookrightarrow \Lambda^\bullet \mathfrak{A}^* \otimes S_{\lambda+1/2+\alpha-\log}(A(t))$ (as the constant functions in degree zero) is a quasi-isomorphism. This is the elementary calculation of the weighted $L^2$ cohomology of a half-line with exponential weights (see Lemma 3.5).

## 1. Preliminaries

This section sets up the necessary notation, defines the two types of cohomology groups and collects some useful facts about $X$ and $\overline{X}$.

*Notational conventions.* Algebraic groups defined over the rational numbers will be denoted by boldface letters ($\mathbf{G}, \mathbf{P}$ etc.), their groups of real points by the corresponding Italic letters ($G, P$ etc.) and the Lie algebras of the latter by Gothic letters ($\mathfrak{G}, \mathfrak{P}$ etc.). Complexes will be denoted $A^\bullet, B^\bullet$ except if they have zero differential when they may be denoted $A^*, B^*$ etc. For example, $H^*(C^\bullet)$ is the cohomology of the complex $C^\bullet$. Complexes of sheaves will be denoted $\mathbf{A}^\bullet, \mathbf{B}^\bullet$ etc. and if $\mathbf{A}^\bullet$ is defined via a particular presheaf that presheaf will generally be denoted $P\mathbf{A}^\bullet$. The sheafification functor is denoted Sh.

1.1. *Notation.* Let $\mathbf{G}, G, K$ be as in the introduction and $D = G/K$. Let $\mathfrak{G} = \text{Lie}(G)$ and $\mathfrak{K} = \text{Lie}(K)$. Let $\mathcal{U}(\mathfrak{G})$ be the universal enveloping algebra of $\mathfrak{G}$ and let $\mathcal{Z}(\mathfrak{G})$ be its centre. For $D \in \mathcal{U}(\mathfrak{G})$ let $R_D$ (resp. $L_D$) be the associated left (resp. right) invariant differential operator on functions on $G$.

Let $\Gamma \subset \mathbf{G}(\mathbb{Q})$ be a neat arithmetic subgroup. (Any arithmetic group has a neat subgroup of finite index ([2, 17.4]).) The locally symmetric space $X = \Gamma \backslash D$ is then a smooth manifold. Let $\text{pr}_\Gamma : D \to X$ and $\text{pr}_K : \Gamma \backslash G \to X$ be the projections.

The group of $\mathbb{Q}$-rational characters of an algebraic group $\mathbf{H}$ defined over $\mathbb{Q}$ is denoted $X^*(\mathbf{H})$. Let $^0\mathbf{H} := \bigcap_{\chi \in X^*(\mathbf{H})} \ker(\chi^2)$. It contains any compact or arithmetic subgroup of $H = \mathbf{H}(\mathbb{R})$. For a subgroup $H \subset G$ let $\Gamma_H := \Gamma \cap H$. For a rational parabolic subgroup $\mathbf{P} \subset \mathbf{G}$ let

$\mathbf{N_P} = $ unipotent radical of $\mathbf{P}$, $N_P = \mathbf{N_P}(\mathbb{R})$ and $\mathfrak{N}_P = \text{Lie}(N_P)$,

$\mathbf{L_P} = \mathbf{P}/\mathbf{N_P}$ the Levi quotient, $L_P = \mathbf{L_P}(\mathbb{R})$,

$\mathbf{S_P} = \mathbb{Q}$-split part of the centre of $\mathbf{L_P}$,

$A_P = \mathbf{S_P}(\mathbb{R})^0$, $\mathfrak{A}_P = \text{Lie}(A_P)$ and $d(P) = \dim(A_P)$,

$\mathbf{M_P} = {}^0\mathbf{L_P}$ and $\Gamma_{M_P}$ is the projection of $\Gamma_P$ to $L_P$.

By a slight abuse of notation $X^*(\mathbf{S_P})$ will be denoted $X^*(A_P)$. We have $\mathfrak{A}_P^* = X^*(A_P) \otimes \mathbb{R}$. Note that $P$ is the semidirect product of $^0P$ and $A_P$. Let $K_P = K \cap P$.

The choice of $x \in D$ (equivalently, of a maximal compact subgroup $K(x) \subset G$) gives canonical lifts $A_P(x), M_P(x) \subset P$ ([6, §1]) and the Langlands decomposition $P = N_P A_P(x) M_P(x)$.

Fix, once and for all, a minimal rational parabolic subgroup $\mathbf{P}_0$ with Lie algebra $\mathfrak{P}_0$. The corresponding groups will be denoted $\mathbf{S}, A, N$ and Lie algebras $\mathfrak{A}, \mathfrak{N}$. Let $\mathfrak{H} \supset \mathfrak{A}$ be a Cartan subalgebra of $\mathfrak{G}_\mathbb{C}$ and $\mathfrak{B} \subset (\mathfrak{P}_0)_\mathbb{C}$ a Borel subalgebra containing $\mathfrak{H}$. Let $\rho$ be half the sum of the positive roots of $\mathfrak{H}$ in $\mathfrak{G}$ in the positive system determined by $\mathfrak{B}$, and let $\rho_P := \rho|\mathfrak{A}_P, \rho_0 := \rho|\mathfrak{A}$.

The choice of $\mathbf{P}_0$ determines a positive system among the roots of $\mathbf{S}$ in $\mathbf{G}$; let $\Delta$ be the set of simple roots. A rational parabolic subgroup $\mathbf{P} \supset \mathbf{P}_0$ corresponds to a subset $\Delta_P \subset \Delta$ (such that $\Delta_P \mapsto \mathbf{P}$ is order-reversing). There is an identification $\mathbf{S}_\mathbf{P} = \left(\bigcap_{\alpha \in \Delta - \Delta_P} \ker(\alpha)\right)^0 \subset \mathbf{S}$. The restrictions $\{\alpha|\mathbf{S}_\mathbf{P} \,|\, \alpha \in \Delta_P\}$ give a basis for $X^*(A_P) \otimes \mathbb{Q}$. The positive cone in $\mathfrak{A}_P^* = X^*(A_P) \otimes \mathbb{R}$ generated by these elements is denoted by $^+\mathfrak{A}_P^*$ (by $^+\mathfrak{A}^*$ for $\mathbf{P} = \mathbf{P}_0$) and its closure is $\overline{^+\mathfrak{A}_P^*}$ ($\overline{^+\mathfrak{A}^*}$ for $\mathbf{P} = \mathbf{P}_0$). For $\lambda \in X^*(A)$ the restriction to $A_P \subset A$ will be denoted $\lambda_P$.

Let $E$ be an irreducible finite-dimensional algebraic representation of $\mathbf{G}$ with highest weight $\Lambda$. There is an admissible inner product ([7, II.2.2]) on $E$ that is unique up to scalars and makes $\mathbb{E}$ a metrized local system.

Given a Lie algebra $\mathfrak{L}$, a reductive group $J$ with $\mathfrak{J} = \mathrm{Lie}(J)$ a subalgebra of $\mathfrak{L}$, and an $(\mathfrak{L}, J)$-module $V$, let $C^\bullet_{(\mathfrak{L}, J)}(V) = \mathrm{Hom}_J(\Lambda^\bullet(\mathfrak{L}/\mathfrak{J}), V)$ be the usual relative Lie algebra cohomology complex ([7, I.1]).

1.2. *Reductive Borel-Serre compactification* ([6], [29, §4], [12, §1]). Fix a maximal compact subgroup $K$ (equivalently, a point $x \in D$) and let $D_P = M_P/K_P$. The equality $P = N_P A_P(x) M_P(x)$ gives a diffeomorphism (depending on $x$) $\mu : N_P \times A_P \times D_P \longrightarrow D$. In these coordinates the action of $A_P$ by $a' \cdot (n, a, z) = (n, a'a, z)$ is the *geodesic action* of Borel-Serre [6]; it is independent of the point $x$ and it commutes with the action of $N_P$, which is given by $n' \cdot (n, a, z) = (n'n, a, z)$ ([3, I, §4]).

The restrictions of the roots in $\Delta_P$ define a diffeomorphism $A_P \to (0, \infty)^{d(P)}$. Then $(0, \infty]^{d(P)}$ is identified with a space $\overline{A}_P$ with an $A_P$-action coming from $A_P \subset \overline{A}_P$. The *corner* associated to $P$ is $D(P) = \overline{A}_P \times_{A_P} D$. Let $e_P := A_P \backslash D$ be the *face* associated to $P$ and $r_P : D \to e_P$ the natural map (the *geodesic retraction*). The face $e_P$ is identified with $^0P/K_P$ and $D(P)$ is identified with $\coprod_{Q \subset P} e_Q$. The union $\widetilde{D} = \coprod_P e_P = \bigcup_P D(P)$ is a space with a properly discontinuous action of $\Gamma$ and $\widetilde{X} := \Gamma \backslash \widetilde{D}$ is the Borel-Serre compactification. The boundary component corresponding to $P$ is $Y_P := \Gamma_P \backslash e_P$. Since $e_P \simeq D_P \times N_P$, $Y_P$ fibres over $X_P := \Gamma_{M_P} \backslash D_P$ with fibres diffeomorphic to the nilmanifold $\Gamma_{N_P} \backslash N_P$. Collapsing these fibres gives the reductive Borel-Serre compactification $\overline{X}$ and a quotient map $\pi : \widetilde{X} \to \overline{X}$. Let $i : X \hookrightarrow \overline{X}$ and $j : X \hookrightarrow \widetilde{X}$ be the inclusions.

Let $\widetilde{e}_P$ be the closure of $e_P$ in $\widetilde{D}$. The retraction $r_P : D \to e_P$ extends to a mapping $r_P : \coprod_{Q \cap P \text{ parabolic}} e_Q \to \widetilde{e}_P$, and the various $r_P$ are compatible. Thus $r_P$ is well-defined on a neighbourhood of the closure $\widetilde{Y}_P$ in $\widetilde{X}$.

1.3. *Two types of neighbourhoods.* Let $x' \in X_P$ be a point in the boundary of $\overline{X}$. Let $O_P$ be a relatively compact set in $D_P$ such that its projection to $X_P$ contains $x'$ and for $t > 0$ let

$$A_P(t) \;=\; \{\, a \in A_P \,|\, a^\alpha > t \text{ for } \alpha \in \Delta_P \,\}.$$

The corresponding subset of $\overline{A}_P$ is $\overline{A}_P(t)$. A fundamental system of neighbourhoods of $x'$ in $\overline{X}$ is given by the *cuspidal neighbourhoods*

$$C_P \;:=\; \mathrm{pr}_\Gamma(\mu(N_P \times \overline{A}_P(t) \times O_P))$$

where $\mu$ is the trivialization given by the choice of $x \in D$ such that $r_P(x)$ projects to $x'$. If there is need to emphasize $t$ and $O_P$ this set will be written $C_P(t, O_P)$.

We shall also need certain neighbourhoods of the closed boundary components in $\widetilde{X}$ that are described in Section 6 of [12]. The following lemma summarizes their properties and gives a clear picture of $\widetilde{X}$ as a manifold-with-corners ([12, 6.4, 6.6]):

LEMMA 1.1. *There exist neighbourhoods $\widetilde{V}_P$ of $\widetilde{Y}_P$ in $\widetilde{X}$ and functions $d_P : \widetilde{V}_P \to [0, \infty]^{d(P)}$ such that (1) $(r_P, d_P) : \widetilde{V}_P \longrightarrow \widetilde{Y}_P \times [0, \infty]^{d(P)}$ is a diffeomorphism, (2) $d_P^{-1}(\infty, \ldots, \infty) = \widetilde{Y}_P$, (3) $d_P \circ \mathrm{pr}_\Gamma$ is $N_P$-invariant and (4) If $P = \cap_{i=1}^{d(P)} Q_i$ for maximal $Q_i$ then $d_P(x) = (d_{Q_1}(x), \ldots, d_{Q_{d(P)}}(x))$.*

Let $V_P = \pi(\widetilde{V}_P)$; it is a neighbourhood of $\overline{X}_P$ in $\overline{X}$. The function $d_P$ makes sense on $V_P$. For $t > 0$ let

$$V_P(t) \;=\; d_P^{-1}((t, \infty]^{d(P)}).$$

The relation of the neighbourhoods $V_P(t)$ to cuspidal neighbourhoods is as follows: For any relatively compact set $O_P \subset X_P$ the intersection

$$V_P(t) \;\cap r_P^{-1}(O_P)$$

is, up to a relatively compact set in $X$, a cuspidal neighbourhood for $P$ and $O_P$.

1.4. *Special differential forms.* A differential form $\sigma \in \Omega^i(X, \mathbb{E})$ is *special* if, for each boundary component $Y_P$ in $\widetilde{X}$, there is a neighbourhood $\widetilde{V}_P(t)$ of $Y_P$ in $\widetilde{X}$ such that on $\widetilde{V}_P(t)$ the form $\sigma$ is lifted (via the retraction $r_P : \widetilde{V}_P \to Y_P$) from a form $\sigma_P$ on $Y_P$ which (when lifted to $e_P$) is $N_P$-invariant. One can make a presheaf on $X$ from this: $\Omega^i_{\mathrm{sp}}(\mathbb{E})$ is the presheaf whose sections over an open set $U \subset X$ are given by: forms $\sigma \in \Omega^i(U \cap X, \mathbb{E})$ such that $\sigma$ is the restriction to $U$ of a special form on $X$. Let $\overline{\Omega}^\bullet_{\mathrm{sp}}(\mathbb{E}) = \mathrm{Sh}(i_*\Omega^\bullet_{\mathrm{sp}}(\mathbb{E}))$ be the associated sheaf on $\overline{X}$.

It is shown in Section 12 of [12] (see also Lemma 2.3 below) that

$$\Omega^\bullet(Y_P, \mathbb{E})^{N_P-\mathrm{inv}} \simeq \Omega^\bullet(X_P, \mathbf{C}^\bullet_{\mathfrak{N}_P}(E)) =: \mathbf{T}^\bullet_P$$

where $\mathbf{C}^\bullet_{\mathfrak{N}_P}(E)$ is the complex of local systems on $X_P$ constructed out of the representation of $L_P = P/N_P$ on the Koszul complex $C^\bullet_{\mathfrak{N}_P}(E) = \Lambda^\bullet \mathfrak{N}_P^* \otimes E$. Here, as in the sequel, the exponent "$N_P$-inv" means that $\mathrm{pr}^*_{\Gamma_P}\sigma$ is left $N_P$-invariant on $e_P$. It is shown in [12, 13.4] that $\overline{\Omega}^\bullet_{\mathrm{sp}}(\mathbb{E})|X_P \cong \mathbf{T}^\bullet_P$.

1.5. *Weighted cohomology complex of sheaves.* The group $A_P$ acts naturally on $C^\bullet_{\mathfrak{N}_P}(E)$. A *weight profile* is the choice, for each $\mathbf{P}$, of an $A_P$-stable subcomplex of $C^\bullet_{\mathfrak{N}_P}(E)$. All profiles occurring in this paper are obtained as follows: Let $p \in X^*(A) \otimes \mathbb{Q}$. Denote the subcomplex of $C^\bullet_{\mathfrak{N}_P}(E)$ on which $A_P$ acts with weights in $p_P + {}^+\mathfrak{A}^*_P$ = open positive cone in $\mathfrak{A}^*_P$ starting at $p_P = p|A_P$ (resp. in $p_P + \overline{{}^+\mathfrak{A}^*_P}$) by $C^\bullet_{\mathfrak{N}_P}(E)_+$ (resp. $C^\bullet_{\mathfrak{N}_P}(E)_-$). These two choices $\mathbf{P} \mapsto C^\bullet_{\mathfrak{N}_P}(E)_+$ and $\mathbf{P} \mapsto C^\bullet_{\mathfrak{N}_P}(E)_-$ are weight profiles and will be denoted $\overline{p}$ ("upper $p$") and $\underline{p}$ ("lower $p$") respectively.

*Remark* 1.2. In [12], a weight profile is defined to be a function $f : \Delta \to \mathbb{Z} + \frac{1}{2}$ and it assigns to $\mathbf{P}$ the subcomplex (denoted $C^\bullet_{\mathfrak{N}_P}(E)_+$ in [12]) with weights $\chi$ such that, for all $\alpha \in \Delta, \chi|A_{Q(\alpha)} > f(\alpha)\alpha$ where $\mathbf{Q}(\alpha)$ is the maximal parabolic subgroup with $\Delta_Q = \{\alpha\}$. Given such a function $f$, let $p := \sum_{\alpha \in \Delta} f(\alpha)\alpha$. Then $\overline{p}$ of the previous paragraph defines the same subcomplex $C^\bullet_{\mathfrak{N}_P}(E)_+$ as $f$ does in [12, §12]. The weight profile $\underline{p}$ corresponds to the function $f - \frac{1}{2}$ (i.e. the $+$-subcomplex defined by $f - \frac{1}{2}$ as in [12, §12] is the subcomplex $C^\bullet_{\mathfrak{N}_P}(E)_-$ above). The notions of weight profile here and in [12] are therefore equivalent.

Denote the local system on $X_P$ associated to $C^\bullet_{\mathfrak{N}_P}(E)_+$ by $\mathbf{C}^\bullet_{\mathfrak{N}_P}(E)_+$. Let $\mathbf{T}^\bullet_{P,+} := \Omega^\bullet(X_P, \mathbf{C}^\bullet_{\mathfrak{N}_P}(E)_+) \subset \mathbf{T}^\bullet_P$. Let $P\mathbf{W}^{\overline{p}}\mathbf{C}^i(\mathbb{E})$ be the presheaf whose sections over $U$ are special forms forms $\sigma$ such that, for each $\mathbf{P}$, there is a $t$ such that $\sigma|U \cap V_P(t)$ is lifted (as in 1.4) from a $\mathbf{C}^\bullet_{\mathfrak{N}_P}(E)_+$-valued form on $X_P$. Then $\mathbf{W}^{\overline{p}}\mathbf{C}^\bullet(\mathbb{E}) := \mathrm{Sh}(P\mathbf{W}^{\overline{p}}\mathbf{C}^\bullet(\mathbb{E}))$ is the weighted cohomology complex of sheaves of profile $\overline{p}$. Its hypercohomology groups $W^{\overline{p}}H^*(\Gamma, E)$ are the weighted cohomology groups of profile $\overline{p}$. The sheaves $\mathbf{W}^{\overline{p}}\mathbf{C}^\bullet(\mathbb{E})$ are fine and so in fact hypercohomology may be replaced by ordinary sheaf cohomology:

$$W^{\overline{p}}H^*(\Gamma, E) = H^*(\mathbf{W}^{\overline{p}}\mathbf{C}^\bullet(\mathbb{E})).$$

The same construction applied to $\underline{p}$ (i.e. starting from $C^\bullet_{\mathfrak{N}_P}(E)_-$) yields a complex of sheaves $\mathbf{W}^{\underline{p}}\mathbf{C}^\bullet(\mathbb{E})$ and its cohomology groups $W^{\underline{p}}H^i(\Gamma, E)$.

*Remark* 1.3. With our conventions, if $p, q \in X^*(A) \otimes \mathbb{Q}$ and $p+q = -2\rho_0$ then the profiles $\overline{p}$ and $\underline{q}$ are dual, so that $\mathbf{W}^{\overline{p}}\mathbf{C}^\bullet(\mathbb{E})$ and $\mathbf{W}^{\underline{q}}\mathbf{C}^\bullet(\mathbb{E}^*)$ are Verdier dual ([12, §§20, 21]) and $W^{\overline{p}}H^i(\Gamma, E) \cong W^{\underline{q}}H^{\dim X - i}(\Gamma, E^*)$.

1.6. *Middle profiles and the Hermitian case.* Let $\mu = \overline{-\rho_0}$ and $\nu = \underline{-\rho_0}$ be the middle profiles of the introduction. If $X$ is Hermitian, the Baily-Borel Satake compactification $\hat{X}$ is a complex projective variety. There is a quotient map $\Phi : \overline{X} \to \hat{X}$ (due to Zucker) that extends the identity on $X$. One of the main results of [12, Th. 23.2] is that both $R\Phi_*\mathbf{W}^\mu\mathbf{C}^\bullet(\mathbb{E})$ and $R\Phi_*\mathbf{W}^\nu\mathbf{C}(\mathbb{E})$ are quasi-isomorphic to the intersection homology complex of sheaves (for the natural stratification) on $\hat{X}$.

1.7. *Weight functions.* A function $w : X \to \mathbb{R}_{>0}$ is called *admissible* if $|R_D w / w|$ is bounded for each $D \in \mathcal{U}(\mathfrak{G})$. A particular family of such functions $w_\lambda$ for $\lambda \in \mathfrak{A}^*$ is described in [11, 2.1] (in the adelic situation). For our purposes, $w_\lambda$ is determined by three properties: (1) it is right $K$-invariant, (2) for $g = nmak$ in a Siegel set for a minimal parabolic subgroup (notation as in the introduction), $w_\lambda(nmak) = O(a^{-\lambda})$ (i.e., both $a^\lambda w_\lambda$ and $a^{-\lambda}(w_\lambda)^{-1}$ are bounded) and (3) if $t$ is large enough then $w_\lambda$ is $N_P$-invariant on $C_P(t, O)$ (i.e., it is $N_P$-invariant near the $P$-boundary). The existence of $w_\lambda$ is an exercise in reduction theory, and while it is certainly not unique, none of the constructions to come are affected by this.

For some fixed dominant weight $\varepsilon$, let $v = \log(w_{-\varepsilon})$. Then $wv$ and $wv^{-1}$ are admissible if $w$ is so.

1.8. *Weighted $L^2$ cohomology complex of sheaves.* In this subsection I shall define a complex of sheaves on $\overline{X}$ which computes (via results of Franke) the weighted $L^2$ cohomology groups. In the sequel I shall deal only with these objects and not the complexes $L^\bullet_{\lambda-\log}(X, \mathbb{E})$ of the introduction.

For $W \subset \Gamma\backslash G$ that is right $K$-stable, let $S_{\lambda-\log}(W)$ (respectively, $S_{\lambda+\log}(W)$) be the space of smooth $K$-finite functions $f : W \to \mathbb{C}$ such that, for any $m \geq 0$ (respectively, for some $m \leq 0$),

$$w_\lambda v^m R_D f \in L^2(W) \text{ for all } D \in \mathcal{U}(\mathfrak{G}).$$

Since $W$ is $K$-stable, $S_{\lambda\pm\log}(W)$ is a $(\mathfrak{G}, K)$-module. Then $U \mapsto S_{\lambda\pm\log}(\mathrm{pr}_K^{-1}(U \cap X)) \otimes E$ defines a presheaf of $(\mathfrak{G}, K)$-modules on $\overline{X}$. Forming the relative Lie algebra complex gives a complex of presheaves by

$$P\mathbf{S}^i_{\lambda\pm\log}(U, \mathbb{E}) = C^i_{(\mathfrak{G}, K)}(S_{\lambda\pm\log}(\mathrm{pr}_K^{-1}(U \cap X)) \otimes E)$$
$$= \mathrm{Hom}_K(\Lambda^i(\mathfrak{G}/\mathfrak{K}), S_{\lambda\pm\log}(\mathrm{pr}_K^{-1}(U \cap X)) \otimes E).$$

Let $\mathbf{S}^i_{\lambda\pm\log}(\mathbb{E}) = \mathrm{Sh}(PS^i_{\lambda\pm\log}(\mathbb{E}))$. The Lie algebra differential $PS^i_{\lambda\pm\log}(U,\mathbb{E}) \to PS^{i+1}_{\lambda\pm\log}(U,\mathbb{E})$ makes $\mathbf{S}^\bullet_{\lambda\pm\log}(\mathbb{E})$ into a complex of sheaves on $\overline{X}$. These sheaves are fine by [29, 4.4]. (Indeed, one can construct partitions of unity for the sheaf $\mathbf{S}^0_{\lambda\pm\log}$, so it is fine. Since the $\mathbf{S}^i_{\lambda\pm\log}$ are modules over $\mathbf{S}^0_{\lambda\pm\log}$ they are also fine.) As a consequence, their hypercohomology is computed by the complex of global sections and so

$$H^*(\mathbf{S}^\bullet_{\lambda\pm\log}(\mathbb{E})) = H^*_{(\mathfrak{G},K)}(S_{\lambda\pm\log}(\Gamma\backslash G)\otimes E).$$

Dropping the powers of $v$ in the $L^2$ condition above gives a subobject $\mathbf{S}^\bullet_\lambda(\mathbb{E}) \subset \mathbf{S}^\bullet_{\lambda\pm\log}(\mathbb{E})$, a subspace $S_\lambda(\Gamma\backslash G) \subset S_{\lambda+\log}(\Gamma\backslash G)$, and groups $H^*_\lambda(X,\mathbb{E})$.

A result of Franke's (of Borel's [4] in the special case of very positive $\lambda$), which follows from Theorem 4 of [11] (applied to the weight functions $w_\lambda v^m$), establishes the relation with the groups $H^i_{\lambda\pm\log}(X,\mathbb{E})$ of the introduction:

PROPOSITION 1.4.   $H^i_{\lambda\pm\log}(X,\mathbb{E}) \cong H^i_{(\mathfrak{G},K)}(S_{\lambda\pm\log}(\Gamma\backslash G)\otimes E).$

The groups $H^*_{\lambda\pm\log}(X,\mathbb{E})$ are always finite-dimensional (because, as we shall see in 3.6, 3.7, the sheaves $\mathbf{S}^\bullet_{\lambda\pm\log}(\mathbb{E})$ are cohomologically constructible on $\overline{X}$), in contrast to the $H^*_\lambda(X,\mathbb{E})$, which need not be so. In case $\lambda = 0$, I shall use the subscript $(2)\pm\log$ instead of $0\pm\log$ (e.g. $\mathbf{S}^\bullet_{(2)\pm\log}(\mathbb{E}), S_{(2)\pm\log}(\Gamma\backslash G), H^*_{(2)\pm\log}(X,\mathbb{E})$ etc.).

*Remark* 1.5.   A Sobolev-type estimate (see [11, Prop. 2.3.2]) implies that $S_{\lambda-\log}(\Gamma\backslash G)$ (resp. $S_{\lambda+\log}(\Gamma\backslash G)$) can also be described as the space of smooth and $K$-finite functions $f$ such that, for any $m \geq 0$ (resp. for some $m \leq 0$) the function $a^{-\lambda-\rho_0}\log(a^\varepsilon)^m R_D f(nmak)$ is bounded on any Siegel set for all $D \in \mathcal{U}(\mathfrak{G})$.

1.9. *Duality.* For $K$-stable $W \subset \Gamma\backslash G$ there is a natural pairing between $S_\lambda(W)$ and $S_{-\lambda}(W)$ given by $(f,g) \mapsto \int_W f\bar{g}$. This also gives a pairing between $S_{\lambda-\log}(W)$ and $S_{-\lambda+\log}(W)$. The pairing is clearly nondegenerate and $G$-invariant, and so (by [7, I.1.5 and I.7.6]) it induces an isomorphism of the cohomology sheaves of $\mathbf{S}^\bullet_{\lambda-\log}(\mathbb{E})$ with those of $\mathbf{S}^\bullet_{-\lambda+\log}(\mathbb{E}^*)$, where $E^*$ is the dual representation of $E$. This gives a duality isomorphism $H^*_{\lambda-\log}(X,\mathbb{E}) \cong H^{\dim X-*}_{-\lambda+\log}(X,\mathbb{E}^*)$.

## 2. Special forms and integrability

The ($N_P$-invariant) sections of $PS^\bullet_{\lambda\pm\log}(\mathbb{E})$ over a cuspidal neighbourhood are studied in 2.1–2.4 and the results are used (2.6) to identify the weighted cohomology sheaves with subsheaves of $\mathbf{S}^\bullet_{\lambda\pm\log}(\mathbb{E})$ for suitable $\lambda$.

2.1. Fix a maximal compact subgroup $K$ for the rest of this section. For a parabolic subgroup $\mathbf{P}$ let $A_P$ and $M_P$ be the canonical lifts determined by $K$. Let $da, dn, dm$ be the left-invariant measures on $A_P, N_P, M_P$ with $dn$ normalized to give volume one to $\Gamma_{N_P} \backslash N_P$. Let $dz$ be the left-$M_P$-invariant volume on $D_P$. The volume form on $D$ is given (up to constants) by $a^{-2\rho_P} dn \wedge da \wedge dz$ and that on $\Gamma \backslash G$ by $a^{-2\rho_P} dn \wedge da \wedge \text{pr}_{K_P}^* dz$.

Let $W \subset \Gamma \backslash G$, let $\omega$ be a smooth $i$-form on $W$, and let $Z \in \Lambda^i(\mathfrak{G}/\mathfrak{K})$. The interior product $i_Z \omega$ of $\omega$ with the left-invariant vector field associated to $Z$ is an $E$-valued function on $W$. Given a smooth $i$-form $\sigma$ on $U \subset X$ the formula

$$\phi_\sigma(Z)(g) = i_Z(\text{pr}_K^* \sigma)(g) \qquad (Z \in \Lambda^i(\mathfrak{G}/\mathfrak{K}),\, g \in G)$$

defines an element $\phi_\sigma \in C^i_{(\mathfrak{G},K)}(C^\infty(\text{pr}_K^{-1}(U)) \otimes E)$. It is well-known (e.g. [7, VII.2.5]) that $\sigma \mapsto \phi_\sigma$ identifies the complex of smooth $\mathbb{E}$-valued forms on $U$ with the complex $C^\bullet_{(\mathfrak{G},K)}(C^\infty(\text{pr}_K^{-1}(U)) \otimes E)$. It follows that, for $U \subset \overline{X}$, the complexes $P\mathbf{S}^\bullet_{\lambda \pm \log}(U, \mathbb{E})$ defined in 1.8 can be considered as subcomplexes of the complex of smooth $\mathbb{E}$-valued forms on $U \cap X$. (This remark is not used explicitly, but it is helpful to keep it in mind.)

2.2. Let $O \subset M_P$ be relatively compact. Define a subset $P(t, O) \subset P$ by

$$P(t, O) := N_P A_P(t) O \subset P.$$

Assume that $t$ is large enough so that $\Gamma$-equivalence and $\Gamma_P$-equivalence are the same on $P(t, O)K$. Let $S_{\lambda-\log}(\Gamma_P \backslash P(t, O) K_P)$ be the space of smooth $K_P$-finite functions $f$ such that, for any $m \geq 0$, $w_\lambda v^m R_D f$ is $L^2$ with respect to the volume $a^{-2\rho_P} dn \wedge da \wedge \text{pr}_{K_P}^* dz$ for every $D \in \mathcal{U}(\mathfrak{P})$. Let $\text{Ind}^{(\mathfrak{G},K)}_{(\mathfrak{P},K_P)}$ be the induction functor taking $(\mathfrak{P}, K_P)$-modules to $(\mathfrak{G}, K)$-modules described in Section 4 of [11] or [7, III.2.1–2.3] (or that of [28, 3.3] at the real place).

LEMMA 2.1.  *There is a canonical identification*

$$S_{\lambda-\log}(\Gamma_P \backslash P(t, O)K) \otimes E = \text{Ind}^{(\mathfrak{G},K)}_{(\mathfrak{P},K_P)} S_{\lambda-\log}(\Gamma_P \backslash P(t, O)K_P) \otimes E.$$

*Proof.* The analogue of the lemma for smooth $K_P$-finite functions on $\Gamma_P \backslash P(t, O)K_P$ is shown in [11, §4]; essentially the same proof applies. The induced module on the right is characterized by a universal property and so it suffices to show that for any $(\mathfrak{G}, K)$-module $V$ with a $(\mathfrak{P}, K_P)$-module map $\varphi : V \to S_{\lambda-\log}(\Gamma_P \backslash P(t, O)K_P)$ there is a $(\mathfrak{G}, K)$-map

$$\tilde{\varphi} : V \to S_{\lambda-\log}(\Gamma_P \backslash P(t, O)K)$$

such that $\varphi$ is the composition of $\tilde{\varphi}$ with the restriction map

$$S_{\lambda-\log}(\Gamma_P \backslash P(t, O)K) \to S_{\lambda-\log}(\Gamma_P \backslash P(t, O)K_P).$$

Define $\tilde{\varphi}$ by $\tilde{\varphi}(v)(pk) = \varphi(k^{-1} \cdot v)(p)$ for $k \in K, p \in P, v \in V$. It is clear that $\tilde{\varphi}(v)$ is smooth and $K$-finite and that $w_\lambda v^m \tilde{\varphi}(v)$ is $L^2$ on $\Gamma_P \backslash P(t,O)K_P$ for every $m \geq 0$. If $D \in \mathfrak{P}$ then $(R_D \tilde{\varphi}(v))(pk) = (R_{\mathrm{Ad}k(D)} \varphi(k^{-1}v))(p)$ and so, for any $m$, $w_\lambda v^m R_D \tilde{\varphi}(v)$ is $L^2$ for all $D \in \mathfrak{P}$. If $D \in \mathfrak{K}$ then $R_D \tilde{\varphi}(v)(pk) = \varphi(D \cdot k^{-1} \cdot v)(p)$ and hence, for any $m$, $w_\lambda v^m R_D \tilde{\varphi}(v)$ is $L^2$ for all $D \in \mathfrak{K}$. Then $\tilde{\varphi}(v) \in S_{\lambda-\log}(\Gamma_P \backslash P(t,O)K)$ showing that $\tilde{\varphi}$ is the desired map and the universal property holds for $S_{\lambda-\log}(\Gamma_P \backslash P(t,O)K)$. $\square$

Given this identification, it follows (see [7, III.2.5]) that:

LEMMA 2.2. *There is a canonical isomorphism of complexes*

$$C^\bullet_{(\mathfrak{G},K)}\big(S_{\lambda-\log}(\Gamma_P \backslash P(t,O)K) \otimes E\big) \;\cong\; C^\bullet_{(\mathfrak{P},K_P)}\big(S_{\lambda-\log}(\Gamma_P \backslash P(t,O)K_P) \otimes E\big).$$

2.3. Suppose that $(A^\bullet, d_A)$ is a complex of $(\mathfrak{L}, J)$-modules. Then by $C^\bullet_{(\mathfrak{L},J)}(A^\bullet)$ we mean the total complex of the double complex with $(i,j)$-term $C^i_{(\mathfrak{L},J)}(A^j)$ and differential given by $d_{(\mathfrak{L},J)} + (-1)^i d_A$ on $C^i_{(\mathfrak{L},J)}(A^j)$.

In the following lemma, $(\mathfrak{L}, J)$ is a pair as in 1.1 and $\mathfrak{L} = \mathfrak{Q} + \mathfrak{R}$ (a sum of Lie algebras) with $\mathfrak{R} \subset \mathfrak{L}$ an ideal. Suppose that $\mathfrak{J} = \mathrm{Lie}(J) \subset \mathfrak{Q}$.

LEMMA 2.3. *Let $V$ be an $(\mathfrak{L}, J)$-module that is a trivial $\mathfrak{R}$-module. There is a canonical isomorphism of complexes*

$$C^\bullet_{(\mathfrak{L},J)}(V \otimes E) \;\cong\; C^\bullet_{(\mathfrak{Q},J)}\big(V \otimes C^\bullet_\mathfrak{R}(E)\big).$$

*Proof.* Since $\Lambda^\bullet(\mathfrak{L}) \cong \Lambda^\bullet \mathfrak{Q} \otimes \Lambda^\bullet \mathfrak{R}$ there is an obvious isomorphism of graded vector spaces between the two complexes and it is easily checked that it carries one differential into the other. $\square$

If $V$ is the $(\mathfrak{P}, K_P)$-module of $N_P$-invariant smooth $K_P$-finite functions on $\Gamma_P \backslash P(t,O)K_P$ then the lemma (with $(\mathfrak{L}, J) = (\mathfrak{P}, K_P)$, $\mathfrak{Q} = \mathfrak{M}_P + \mathfrak{A}_P$ and $\mathfrak{R} = \mathfrak{N}_P$) identifies the smooth $N_P$-invariant $\mathbb{E}$-valued forms on $\mathrm{pr}_{K_P}(\Gamma_P \backslash P(t,O)K_P)$ with the smooth $\mathbf{C}^\bullet_{\mathfrak{N}_P}(E)$-valued forms on the same neighbourhood ([12, §12]). We shall use it for $V = S_{\lambda-\log}(\Gamma_P \backslash P(t,O)K_P)^{N_P-\mathrm{inv}}$.

2.4. *Invariant sections over cuspidal neighbourhoods.* To see how $P\mathbf{S}^\bullet_{\lambda-\log}(\mathbb{E})$ behaves with respect to the product structure of cuspidal neighbourhoods we restrict attention to $N_P$-invariants. For $\delta \in \mathfrak{A}^*_P$, let $S_{\delta-\log}(A_P(t))$ be the space of smooth functions $f$ such that, for any $m \geq 0$,

$$\log(a^\varepsilon)^m a^{-\delta} R_D f \in L^2(A_P(t)) \quad \text{for all } D \in \mathcal{U}(\mathfrak{A}_P)$$

($\varepsilon$ a dominant weight as in 1.7).

LEMMA 2.4. *Under the identification $A_P(t) \times \Gamma_{M_P}\backslash OK_P = \Gamma_{M_P}\backslash A_P(t)OK_P$ there is an inclusion*

$$S_{\lambda_P+\rho_P-\log}(A_P(t)) \otimes S_{(2)}(\Gamma_{M_P}\backslash OK_P) \subset S_{\lambda-\log}(\Gamma_P\backslash P(t,O)K_P)^{N_P-\text{inv}}.$$

*Proof.* The right-hand side is identified with the space of functions $f$ on $\Gamma_P N_P\backslash N_P A_P(t)OK_P$ with $w_\lambda v^m R_D f \in L^2$ for all $D \in \mathcal{U}(\mathfrak{P})$ and all $m \geq 0$. For any $D \in \mathcal{U}(\mathfrak{P})$ the operator $R_D$ can be written as a finite sum of terms like $g(n,l)L_{D_1}R_{D_2}$ where $D_1 \in \mathcal{U}(\mathfrak{N}_P), D_2 \in \mathcal{U}(\mathfrak{L}_P)$ (see e.g. [10, 2.5]). Since $R_{D_2}f$ is $N_P$-invariant if $f$ is $N_P$-invariant and $L_{D_1}f = 0$ for such $f$, the space on the right is identified with the space $S_{\lambda-\log}(\Gamma_{M_P}\backslash A_P(t)OK_P)$ of functions $f$ such that, for all $m$, $a^{-\lambda_P}\log(a^\varepsilon)^m R_D f \in L^2$ for all $D \in \mathcal{U}(\mathfrak{L}_P)$. Since $\mathfrak{L}_P = \mathfrak{A}_P + \mathfrak{M}_P$ and $[\mathfrak{A}_P, \mathfrak{M}_P] = 0$, the expression for the volume form implies the lemma. □

*Remark* 2.5. The inclusion of the lemma extends to Grothendieck's topological tensor product ([7, IX.6.1] or [26]) on the left, and the extension can be shown to be an isomorphism. However, this is not necessary in the sequel.

Let $C_P$ be the cuspidal neighbourhood in $\overline{X}$ associated to $t$ and $\text{pr}_{K_P}(O)$ (the basepoint is given by $K$); then $\text{pr}_K^{-1}(C_P \cap X) = \Gamma_P\backslash P(t,O)K$. Combining Lemmas 2.2 and 2.3 we have canonical isomorphisms

$$P\mathbf{S}^\bullet_{\lambda-\log}(C_P, \mathbb{E})^{N_P-\text{inv}} \cong C^\bullet_{(\mathfrak{P},K_P)}\big(S_{\lambda-\log}(\Gamma_P\backslash NA_P(t)OK_P)^{N_P-\text{inv}} \otimes E\big)$$

$$\cong C^\bullet_{(\mathfrak{M}_P+\mathfrak{A}_P,K_P)}\Big(S_{\lambda-\log}(\Gamma_{M_P}\backslash A_P(t)OK_P) \otimes C^\bullet_{\mathfrak{N}_P}(E)\Big).$$

Let $\mathbb{C}_\alpha$ denote the $A_P$-module with character $\alpha$ and let $V^\bullet_\alpha$ be the complex of $M_P$-modules $C^\bullet_{\mathfrak{N}_P}(E)_\alpha$ considered as trivial $A_P$-modules. Then $C^\bullet_{\mathfrak{N}_P}(E)_\alpha \cong \mathbb{C}_\alpha \otimes V^\bullet_\alpha$ as $(\mathfrak{M}_P + \mathfrak{A}_P, K_P)$-modules. Associated to $\mathbb{C}_\alpha$ in the usual way is a metrized bundle on $A_P(t)$ and hence one on $\Gamma_{M_P}\backslash A_P(t)OK_P$. The $(\mathfrak{M}_P + \mathfrak{A}_P, K_P)$-module $S_{\lambda-\log}(\Gamma_{M_P}\backslash A_P(t)OK_P) \otimes \mathbb{C}_\alpha$ can be identified with the smooth sections of this bundle which satisfy the appropriate integrability condition on all derivatives. Choice of a nowhere-zero section (the bundle is trivial) gives an isomorphism of $(\mathfrak{M}_P + \mathfrak{A}_P, K_P)$-modules

$$S_{\lambda-\log}(\Gamma_{M_P}\backslash A_P(t)OK_P) \otimes \mathbb{C}_\alpha \cong S_{\lambda+\alpha-\log}(\Gamma_{M_P}\backslash A_P(t)OK_P).$$

Decomposing by $A_P$-weights and using this isomorphism we have:

$$P\mathbf{S}^\bullet_{\lambda-\log}(C_P, \mathbb{E})^{N_P-\text{inv}} \cong \bigoplus_\alpha C^\bullet_{(\mathfrak{M}_P+\mathfrak{A}_P,K_P)}\big(S_{\lambda+\alpha-\log}(\Gamma_{M_P}\backslash A_P(t)OK_P)) \otimes V^\bullet_\alpha\big).$$

This is in fact a decomposition as complexes because the differential on $C^\bullet_{\mathfrak{N}_P}(E)$ commutes with the $A_P$-action. Now by Lemma 2.3 (with $(\mathfrak{L}, J) = (\mathfrak{M}_P + \mathfrak{A}_P, K_P), \mathfrak{R} = \mathfrak{A}_P, \mathfrak{Q} = \mathfrak{M}_P$) this is canonically isomorphic to

$$\bigoplus_\alpha C^\bullet_{(\mathfrak{M}_P,K_P)}\Big(C^\bullet_{\mathfrak{A}_P}(S_{\lambda+\alpha-\log}(\Gamma_{M_P}\backslash A_P(t)OK_P)) \otimes V^\bullet_\alpha\Big).$$

2.5. *A subcomplex.* Let $C^\infty_{\text{res}}(\Gamma_{M_P}\backslash OK_P)$ be the smooth $K_P$-finite functions on $\Gamma_{M_P}\backslash OK_P$ that are restrictions of smooth $K_P$-finite functions on $\Gamma_{M_P}\backslash M_P$. (Note that such functions belong to $S_{(2)}(\Gamma_{M_P}\backslash OK_P)$.) Let $A_\alpha$ be the constant functions on $A_P(t)$ if $\alpha \in p_P + {}^+\mathfrak{A}^*_P$ and 0 otherwise. Suppose that $\lambda_P \in -p_P - \rho_P + \overline{{}^+\mathfrak{A}^*_P}$. Then $A_\alpha \subset S_{\lambda_P + \rho_P + \alpha - \log}(A_P(t))$ and so

$$A_\alpha \otimes C^\infty_{\text{res}}(\Gamma_{M_P}\backslash OK_P) \ \subset \ S_{\lambda_P + \rho_P + \alpha - \log}(A_P(t)) \otimes S_{(2)}(\Gamma_{M_P}\backslash OK_P).$$

By Lemma 2.4 there is an inclusion (in degree zero)

$$A_\alpha \otimes C^\infty_{\text{res}}(\Gamma_{M_P}\backslash OK_P) \ \subset \ C^\bullet_{\mathfrak{A}_P}\bigl(S_{\lambda + \alpha - \log}(\Gamma_{M_P}\backslash OK_P)\bigr).$$

Then the complex

$$W^\bullet_p(C_P, \mathbb{E}) \ = \ \bigoplus_\alpha C^\bullet_{(\mathfrak{M}_P, K_P)}\bigl(A_\alpha \otimes C^\infty_{\text{res}}(\Gamma_{M_P}\backslash OK_P) \otimes V^\bullet_\alpha\bigr)$$

is a subcomplex of $P\mathbf{S}^\bullet_{\lambda - \log}(C_P, \mathbb{E})^{N_P - \text{inv}}$. (Note that the second factor in the $\alpha$-summand is identified in the standard way with the smooth $\mathbb{V}^\bullet_\alpha$-valued forms on $\Gamma_{M_P}\backslash \text{pr}_{K_P}(O) \subset X_P$ which are restrictions of $\mathbb{V}^\bullet_\alpha$-valued forms on $X_P$.)

2.6. *The inclusion.* Let $V_P(t)$ be as in 1.2. We write $P \succ R$ if $P$ contains a $\Gamma$-conjugate of $R$ (equivalently, if $\overline{X}_P \supset X_R$). For each $P$, let

$$V'_P(t) = V_P(t) - \bigcup_{R \prec P} V_R(t).$$

Then $V_R(t)$ is the disjoint union of $V'_P(t)$ for $P \succ R$ and $r_P(V'_P(t)) \subset X_P$ is relatively compact.

PROPOSITION 2.6. *If $\lambda \in \mathfrak{A}^*$ satisfies $\lambda + \rho_0 + p \in \overline{{}^+\mathfrak{A}^*}$ then $\mathbf{W}^p\mathbf{C}^\bullet(\mathbb{E})$ is identified with a subsheaf of $\mathbf{S}^\bullet_{\lambda - \log}(\mathbb{E})$.*

*Proof.* Let $U \subset \overline{X}$ be open and let $\sigma \in P\mathbf{W}^{\overline{p}}\mathbf{C}^i(U, \mathbb{E})$. By definition $\sigma$ is the restriction to $U$ of a special form (also denoted $\sigma$) on all of $X$. Let $U'_{P,t} = U \cap V'_P(t)$. Then $U$ is the union of the $U'_{P,t}$ and a relatively compact set $U_{G,t}$ in $X$. Choose $O_P \subset M_P$ projecting to an $\varepsilon$-neighbourhood of $r_P(U'_{P,t})$ in $X_P$ and let $C_P$ be the associated cuspidal neighbourhood. By increasing $t$ we may assume that $\sigma|C_P \cap X$ is lifted from an $N_P$-invariant form $\sigma_P$ (with profile $p$) on $C_P \cap Y_P$. Then on $C_P \cap X$, $\phi_\sigma$ (notation from 2.1) belongs to $W^i_p(C_P, \mathbb{E})$, which is a subcomplex of $P\mathbf{S}^i_{\lambda - \log}(C_P, \mathbb{E})$ provided $\lambda_P + \rho_P + \alpha \in \overline{{}^+\mathfrak{A}^*_P}$. The various $C_P$ together with a relatively compact set in $X$ form a finite cover of $U$ and $\phi_\sigma$ is a section of $P\mathbf{S}^i_{\lambda - \log}(\mathbb{E})$ over each $C_P$. Then $\phi_\sigma \in P\mathbf{S}^i_{\lambda - \log}(U, \mathbb{E})$. We have shown $P\mathbf{W}^{\overline{p}}\mathbf{C}^\bullet(\mathbb{E}) \subset P\mathbf{S}^\bullet_{\lambda - \log}(\mathbb{E})$ and the proposition follows. □

We have actually identified the map $\mathbf{W}^{\overline{p}}\mathbf{C}^\bullet(\mathbb{E}) \to \mathbf{S}^\bullet_{\lambda - \log}(\mathbb{E})$ on the presheaf level; this will be useful in the sequel. The arguments of 2.1–2.5, with the obvious changes, give an injection $\mathbf{W}^{\underline{p}}\mathbf{C}^\bullet(\mathbb{E}) \to \mathbf{S}^\bullet_{\lambda + \log}(\mathbb{E})$ for the same $\lambda$.

## 3. Quasi-isomorphism

In 3.1–3.6 various results are established which are put together in 3.7 to show that the inclusion of the previous section induces an isomorphism in stalk cohomology at any point for $\lambda$ as in Theorem A.

3.1. *The constant term.* Let $S \subset D$ be an $N_P$-stable set. Given a differential form $\sigma$ on $\Gamma_P \backslash S$ (or on $\Gamma \backslash S$) we can take the *constant term of $\sigma$ along $P$*:

$$CT_P(\sigma)(x) := \int_{\Gamma_{N_P} \backslash N_P} \sigma(nx) dn$$

where $dn$ is the normalized Haar measure on $N_P$ as before. The constant term of $\sigma$ is an $N_P$-invariant (i.e., $\mathrm{pr}^*_{\Gamma_P} \sigma$ is left $N_P$-invariant) differential form on $\Gamma_P \backslash S$ or a form on $\Gamma_P N_P \backslash S$. The same formula defines the constant term of a function $f$ on $\Gamma \backslash G$; it is an $N_P$-invariant function on $\Gamma_P \backslash G$ or a function on $\Gamma_P N_P \backslash G$.

3.2. *Reduction to $N_P$-invariants.* Let **P** be a rational parabolic subgroup and $O \subset M_P$ a relatively compact subset. Let $P(t_0, O) \subset P$ be defined as in 2.2. Let

$$B(\Gamma_P \backslash P(t_0, O)K) \subset S_{\lambda-\log}(\Gamma_P \backslash P(t_0, O)K)$$

be the subspace of functions $f$ such that $CT_P(f)$ vanishes on $\Gamma_P \backslash P(t, O)K$ for some $t > t_0$. This condition on the constant term is preserved by $\mathcal{U}(\mathfrak{G})$ and $K$ since they act on the right and so $B(\Gamma_P \backslash P(t_0, O)K)$ is a $(\mathfrak{G}, K)$-module. It is the kernel of the map of $(\mathfrak{G}, K)$-modules

$$S_{\lambda-\log}(\Gamma_P \backslash P(t_0, O)K) \longrightarrow \varinjlim_t S_{\lambda-\log}(\Gamma_P \backslash P(t, O)K)^{N_P-\mathrm{inv}}$$

given by the constant term.

LEMMA 3.1. *This map is surjective.*

*Proof.* An element of the right-hand side is a pair $(f_P, t)$ with $f_P \in S_{\lambda-\log}(\Gamma_P \backslash P(t_0, O)K)^{N_P-\mathrm{inv}}$. Let $\varphi_P$ be a smooth function on $A_P(t_0)$ such that $\varphi_P \equiv 1$ on $A_P(t+1)$ and $\varphi \equiv 0$ on the complement of $A_P(t)$. Define $\varphi'_P$ on $\Gamma_P \backslash P(t_0, O)K$ by $\varphi'_P(nmak) = \varphi_P(a)$. Then $\varphi'_P \equiv 1$ on $\Gamma_P \backslash P(t+1, O)K$, $\varphi'_P \equiv 0$ on the complement of $\Gamma_P \backslash P(t, O)K$, and $\varphi'_P$ is left $N_P$-invariant and right $K$-invariant. It belongs to $S_{(2)}(\Gamma_P \backslash P(t_0, O)K)$ by Lemma 2.4. Now let $f = f_P \varphi'_P$ on $P(t_0, O)$. Then $CT_P(f) = f_P$. □

Fix $x \in X_P$. By Lemma 3.1 there is a surjection of $(\mathfrak{G}, K)$-modules

$$\varinjlim_{t,O} S_{\lambda-\log}(\Gamma_P \backslash P(t, O)K) \longrightarrow \varinjlim_{t,O} S_{\lambda-\log}(\Gamma_P \backslash P(t, O)K)^{N_P-\mathrm{inv}}$$

where $t \to \infty$ and $\text{pr}_{K_P}(O)$ runs through a fundamental system of neighbourhoods of $x$ in $X_P$. It has kernel $\varinjlim_{t,O} B(\Gamma_P \backslash P(t,O)K)$. This gives a surjection (which we denote CT) of complexes

$$C^\bullet_{(\mathfrak{G},K)}\big(\varinjlim_{t,O} S_{\lambda-\log}(\Gamma_P\backslash P(t,O)K) \otimes E\big) \xrightarrow{\text{CT}}$$
$$C^\bullet_{(\mathfrak{G},K)}\big(\varinjlim_{t,O} S_{\lambda-\log}(\Gamma_P\backslash P(t,O)K)^{N_P-\text{inv}} \otimes E\big)$$

in which the first term is exactly the stalk of $\mathbf{S}^\bullet_{\lambda-\log}(\mathbb{E})$ at $x$. CT has kernel $C^\bullet_{(\mathfrak{G},K)}\big(\varinjlim_{t,O} B(\Gamma_P\backslash P(t,O)K) \otimes E\big)$.

PROPOSITION 3.2.  CT *is a quasi-isomorphism.*

*Proof.* The argument for this in [8, 1.2] for $\lambda = 0$ (see [19, Lemma 2.4] for SL(2)) works for any $\lambda$. For the sake of completeness, I shall repeat it here. By a Shapiro lemma ([7, III.2.5]) the kernel of CT has cohomology

$$H^*_{(\mathfrak{G},K)}\big(\varinjlim_{t,O} B(\Gamma_P\backslash P(t,O)K) \otimes E\big) \cong H^*_{(\mathfrak{P},K_P)}\big(\varinjlim_{t,O} B(\Gamma_P\backslash P(t,O)K_P) \otimes E\big).$$

The Hochschild-Serre sequence converging to this has $E_2$ term

$$H^*_{(\mathfrak{M}_P+\mathfrak{A}_P,K_P)}\big(H^*_{\mathfrak{N}_P}\big(\varinjlim_{t,O} B(\Gamma_P\backslash P(t,O)K_P) \otimes E\big)\big).$$

It suffices to show that the $\mathfrak{N}_P$-cohomology vanishes. Let $N_P = N_0 \supset N_1 \supset \cdots \supset N_{n+1} = \{e\}$ be the descending central series of $N_P$. Let $\Gamma_r = \Gamma \cap N_r$ and let $\mathfrak{N}_r = \text{Lie}(N_r) \subset \mathfrak{N}_P$. Averaging functions over $\Gamma_r \backslash N_r$ defines a map

$$\text{CT}_r : \varinjlim_{t,O} S_{\lambda-\log}(\Gamma_P\backslash P(t,O)K_P) \to \varinjlim_{t,O} S_{\lambda-\log}(\Gamma_P\backslash P(t,O)K_P)^{N_r-\text{inv}}.$$

Then $\ker(\text{CT}_0) = \varinjlim_{t,O} B(\Gamma_P\backslash P(t,O)K_P) \supset \ker(\text{CT}_1) \supset \cdots \supset \ker(\text{CT}_{n+1}) = 0$. Let $B_r^{r+1} := \text{im}(\text{CT}_{r+1}) \cap \ker(\text{CT}_r)$. Then $\ker(\text{CT}_r)/B_r^{r+1} \cong \ker(\text{CT}_{r+1})$ for each $r$. Suppose we know that

$$(*) \qquad\qquad H^*_{\mathfrak{N}_r/\mathfrak{N}_{r+1}}(B_r^{r+1} \otimes H^*_{\mathfrak{N}_{r+1}}(E)) = 0$$

for each $r$. In particular, for $r = n$, $H^*_{\mathfrak{N}_n}(\ker(\text{CT}_n) \otimes E) = 0$. Assume that $H^*_{\mathfrak{N}_r}(\ker(\text{CT}_r) \otimes E) = 0$ for some $r < n$. By the Hochschild-Serre sequence, it follows that $H^*_{\mathfrak{N}_{r-1}}(\ker(\text{CT}_r) \otimes E) = 0$. Since $\ker(\text{CT}_r)/B_r^{r+1} \cong \ker(\text{CT}_{r+1})$, $(*)$ implies that $H^*_{\mathfrak{N}_{r-1}}(\ker(\text{CT}_{r-1}) \otimes E) = 0$. By downward induction on $r$, we have shown that $H^*_{\mathfrak{N}_0}(\ker(\text{CT}_0) \otimes E) = 0$, proving the proposition.

It remains to prove $(*)$. Since $\mathfrak{N}_r$ acts nilpotently on $H^*_{\mathfrak{N}_{r+1}}(E)$, it will suffice to show $H^*_{\mathfrak{N}_r/\mathfrak{N}_{r+1}}(B_r^{r+1}) = 0$. (Indeed, there is a filtration of $H^*_{\mathfrak{N}_{r+1}}(E)$ by $\mathfrak{N}_r/\mathfrak{N}_{r+1}$-submodules so that $\mathfrak{N}_r/\mathfrak{N}_{r+1}$ acts trivially on the associated graded module. This gives a filtration of $B_r^{r+1} \otimes H^*_{\mathfrak{N}_{r+1}}(E)$ by $\mathfrak{N}_r/\mathfrak{N}_{r+1}$-submodules for which the associated graded module is a sum of copies of $B_r^{r+1}$, whence the assertion.) To prove this it is standard to use a Fourier series argument. Expanding functions in $B_r^{r+1}$ in Fourier series on the compact torus $T = N_r/N_{r+1}\Gamma_r$ (with Lie algebra $\mathfrak{T} = \mathrm{Lie}(T) = \mathfrak{N}_r/\mathfrak{N}_{r+1}$), we write

$$f(\exp(Z)\,g) = \sum_{\xi \in L - \{0\}} f_\xi(g)\exp(\xi(Z)) \qquad (Z \in \mathfrak{T}, g \in P(t,O)K_P)$$

for some lattice $L \subset \mathfrak{T}^*$. The term $\xi = 0$ does not appear because $f \in \ker(\mathrm{CT}_r)$. Choose an inner product on $\mathfrak{T}$ and a basis $\{Z_i\}$ orthonormal with respect to it. For $Z \in \mathfrak{T} = \mathfrak{N}_r/\mathfrak{N}_{r+1}$ let $R_Z$ denote the differential operator on smooth $K_P$-finite left $N_{r+1}$-invariant functions on $P(t,O)K_P$. Then $(R_{Z_i}f_\xi)(g) = \xi(\mathrm{Ad}(g)Z_i)f_\xi(g)$. Set $\delta(\xi,g) = \sum_i \xi(\mathrm{Ad}(g)Z_i)^2$. As $g$ varies in $P(t,O)K_P$ and $\xi \in L - \{0\}$, the function $\delta(\xi,g)^{-1}$ and all its derivatives with respect to $\mathcal{U}(\mathfrak{P})$ are smooth and bounded on $P(t,O)K_P$. Following [8, 1.2], define a Green's operator $G: B_r^{r+1} \to B_r^{r+1}$ by

$$Gf(g) = \sum_\xi f_\xi(g)\,\delta(\xi,g)^{-1}.$$

Now $C = -\sum_i R_{Z_i}^2$ is the Casimir operator on $B_r^{r+1}$. It is clear from the definition that $G$ is bounded and $GC = \mathrm{id}$. It is well-known that if $C$ has a bounded inverse on an $\mathfrak{T}$-module $V$ of smooth vectors then $H_\mathfrak{T}^*(V)$ must vanish. (Indeed, for $d$ the differential on $C_\mathfrak{T}^\bullet(V)$, $d^*$ its adjoint, and $\Delta = dd^* + d^*d$, the invertibility of $C$ implies that of $\Delta$ (by Kuga's formula) and then $\Delta^{-1}d^*: C_\mathfrak{T}^\bullet(V) \to C_\mathfrak{T}^{\bullet-1}(V)$ is a bounded homotopy of $d$ with id.) $\square$

*Remark* 3.3. Analogous assertions are in [9, 4.7], [29, 4.24]; the original result of this kind is due to van Est [27].

3.3. *Stalk of weighted cohomology.* Let $C_P(t,O)$ be the cuspidal neighbourhood associated to $t$ and $\mathrm{pr}_{K_P}(O) \subset D_P$. (The dependence on $t,O$ is emphasized as I want to vary them now.) Let $W_p^\bullet(C_P(t,O), \mathbb{E})$ be the complex of lifted forms on $C_P(t,O)$ as defined in 2.2. If $\sigma \in P\mathbf{W}^{\overline{p}}\mathbf{C}^i(C_P(t_0,O), \mathbb{E})$ then its constant term along $P$ belongs to $W_p^i(C_P(t,O), \mathbb{E})$ for some $t$. This defines a map

$$P\mathbf{W}^{\overline{p}}\mathbf{C}^\bullet(C_P(t_0,O), \mathbb{E}) \to \varinjlim_t W_p^\bullet(C_P(t,O), \mathbb{E}).$$

Varying $t,O$ as before gives an identification of the stalk of $\mathbf{W}^{\overline{p}}\mathbf{C}^\bullet(\mathbb{E})$ at $x$:

$$\varinjlim_{t,O} P\mathbf{W}^{\overline{p}}\mathbf{C}^\bullet(C_P(t,O), \mathbb{E}) \;\cong\; \varinjlim_{t,O} W_p^\bullet(C_P(t,O), \mathbb{E}).$$

3.4. In the following lemma $(\mathfrak{L}, J)$ is as in 1.1 and $J$ is assumed compact.

LEMMA 3.4. *Let $A^\bullet$ be a complex of $(\mathfrak{L}, J)$-modules and $A_0^\bullet \subset A^\bullet$ a subcomplex of $(\mathfrak{L}, J)$-submodules such that the inclusion $A_0^\bullet \to A^\bullet$ induces an isomorphism $H^*(A_0^\bullet) \to H^*(A^\bullet)$ of $(\mathfrak{L}, J)$-modules. Then $C_{(\mathfrak{L},J)}^\bullet(A_0^\bullet) \to C_{(\mathfrak{L},J)}^\bullet(A^\bullet)$ is a quasi-isomorphism.*

*Proof.* Let $T^\bullet$ be the total complex of $C^\bullet_{(\mathfrak{L},J)}(A^\bullet)$ and $T_0^\bullet$ the subcomplex corresponding to $A_0^\bullet$. Define a filtration of $T^\bullet$ by letting $F^i(T^\bullet)$ be the total complex of the double subcomplex with $p,q$-term equal to $C^p_{(\mathfrak{L},J)}(A^q)$ if $p \leq n$ and zero if $p > n$. The spectral sequence of this filtration has $E_1$-term

$$E_1^{p,q} = H^q(C^p_{(\mathfrak{L},J)}(A^\bullet)).$$

$F^i$ induces a filtration of $T_0^\bullet$; the $E_1$-term of the corresponding spectral sequence is $H^q(C^p_{(\mathfrak{L},J)}(A_0^\bullet))$. Then $A_0^\bullet \to A^\bullet$ induces

$$H^q(C^p_{(\mathfrak{L},J)}(A_0^\bullet)) \longrightarrow H^q(C^p_{(\mathfrak{L},J)}(A^\bullet))$$

on the $E_1$-terms. Now $C^p_{(\mathfrak{L},J)}(A^\bullet) = (\Lambda^p(\mathfrak{L}/\mathfrak{J})^* \otimes A^\bullet)^J$. Since tensoring with $\Lambda^p(\mathfrak{L}/\mathfrak{J})^*$ is exact and taking $J$-invariants is exact ($J$ is compact),

$$H^q(C^p_{(\mathfrak{L},J)}(A^\bullet)) \cong C^p_{(\mathfrak{L},J)}(H^q(A^\bullet))$$

and similarly for $A_0^\bullet$. The $E_1$-differential becomes $d_{(\mathfrak{L},J)} : C^p_{(\mathfrak{L},J)}(H^q(A^\bullet)) \to C^{p+1}_{(\mathfrak{L},J)}(H^q(A^\bullet))$ under this isomorphism. The inclusion $A_0^\bullet \to A^\bullet$ induces $C^p_{(\mathfrak{L},J)}(H^q(A_0^\bullet)) \cong C^p_{(\mathfrak{L},J)}(H^q(A^\bullet))$, which is, by hypothesis, an isomorphism on the $E_1$-terms. This implies the lemma. $\square$

3.5. *A crucial lemma.* Let $A'_\alpha$ be the space of constant functions on $A_P(t)$ if $\lambda_P + \rho_P + \alpha \in {}^+\mathfrak{A}_P^*$ and zero otherwise. As explained in 2.4, there is an isomorphism of $(\mathfrak{A}_P + \mathfrak{M}_P, K_P)$-modules

$$S_{\lambda-\log}(\Gamma_{M_P}\backslash A_P(t)OK_P) \otimes \mathbb{C}_\alpha \;\cong\; S_{\lambda+\alpha-\log}(\Gamma_{M_P}\backslash A_P(t)OK_P).$$

There is also an inclusion of $(\mathfrak{A}_P + \mathfrak{M}_P, K_P)$-modules

$$A'_\alpha \otimes S_{(2)}(\Gamma_{M_P}\backslash OK_P) \;\subset\; S_{\lambda+\alpha-\log}(\Gamma_{M_P}\backslash A_P(t)OK_P).$$

LEMMA 3.5. *The inclusion (in degree zero)*

$$A'_\alpha \otimes S_{(2)}(\Gamma_{M_P}\backslash OK_P) \;\subset\; C^\bullet_{\mathfrak{A}_P}(S_{\lambda+\alpha-\log}(\Gamma_{M_P}\backslash A_P(t)OK_P))$$

*induces an isomorphism of $(\mathfrak{M}_P, K_P)$-modules in cohomology.*

*Proof.* Let $\alpha_1, \ldots, \alpha_d$ be the restrictions to $A_P$ of the simple roots that are nontrivial on $A_P$. Let $\check{\alpha}_i$ be the dual basis of $\mathfrak{A}_P$. Then $\mathfrak{A}_P = \sum_i \mathfrak{A}_i$ for $\mathfrak{A}_i = \mathbb{R}\check{\alpha}_i$. Set $\mathfrak{A}_{<i} = \sum_{j<i} \mathfrak{A}_j$. Let $A_i, A_{<i} \subset A_P$ be the subgroups given by $\mathfrak{A}_i, \mathfrak{A}_{<i}$ respectively so that $A_P = \times_i A_i$. Let $r_i = \langle \rho_P + \lambda_P + \alpha, \check{\alpha}_i\rangle$. Let us first assume that $\rho_P + \lambda_P + \alpha \in {}^+\mathfrak{A}_P^*$, i.e., all $r_i > 0$.

Let $S := S_{\lambda+\alpha-\log}(\Gamma_{M_P}\backslash A_P(t)OK_P)$ for the rest of the proof. Let $S_d \subset S$ be the space of functions that are $A_d$-constant. Consider the inclusion $S_d \subset \Lambda^\bullet \mathfrak{A}_d^* \otimes S$ in degree zero. The latter is the two-term complex

$$\Lambda^0 \mathfrak{A}_d^* \otimes S \;\to\; \Lambda^1 \mathfrak{A}_d^* \otimes S$$

with differential $f \mapsto \alpha_d \otimes R_d f$, where $R_d$ is the left-invariant differential operator associated to $\check{\alpha}_d \in \mathfrak{A}_P$. Note that since $A_P$ commutes with $K_P$ and with $O$ we have $R_d f = L_d f$. For $f \in S$ let $If$ be defined by

$$If(a_{<d}, a_d, o, k) = -\int_{a_d}^{\infty} f(a_{<d}, u, o, k)\, du$$

for $a_{<d} \in A_{<d}, a_d \in A_d, o \in O, k \in K_P$, and $du$ the Haar measure on $A_d$. Then $R_d If = L_d If = f$. So $\alpha_d \otimes f \mapsto If$ gives an inverse to the differential of $\Lambda^\bullet \mathfrak{A}_d^* \otimes S$ provided $If \in S$. Let us see that this is so. Let $s$ be defined by $a_d = \exp(s.\check{\alpha}_d)$ and write $f(a_{<d}, s, o, k)$ for $f(a_{<d}, \exp(s.\check{\alpha}_d), o, k)$. It suffices to show that, for each $a_{<d}, o, k$ and $m \geq 0$ the function $e^{-r_i s} s^m R_d^l If(a_{<d}, s, o, k)$ is $L^2$ (with respect to the measure $ds$) on $(\log(t), \infty)$ for all nonnegative $l$. Note that $R_d$ is simply the derivative with respect to $s$. Then $f \in S$ implies that for each $a_{<d}, o, k$ and $m \geq 0$ the function $e^{-r_i s} s^m f^{(l)}(a_{<d}, s, o, k)$ is $L^2$ on $(\log(t), \infty)$ for all $l \geq 0$. By the Sobolev inequality it follows that, for any $m \geq 0$, $e^{-r_i s} s^m f^{(l)}$ is bounded (in $s$) for all $l \geq 0$. But then, for each $m$, the function $e^{-r_i s} R_d^l f$ is bounded by $s^{-m}$ for all $l$. Hence for each $a_{<d}, o, k, m$, the function $e^{-r_i s} s^m R_d^l If(a_{<d}, s, o, k)$ is $L^2$ for all $l \geq 0$. Then $If \in S$. As remarked earlier, $d(If) = \alpha_d \otimes f$ and so the complex $\Lambda^\bullet \mathfrak{A}_d^* \otimes S$ has zero cokernel. It clearly has kernel equal to $\Lambda^0 \mathfrak{A}_d^* \otimes S_d$. Therefore the induced map in cohomology $S_d \to H^*_{\mathfrak{A}_d}(S)$ is an isomorphism; furthermore, it is an isomorphism of $\mathfrak{A}_{<d} = \mathfrak{A}_P/\mathfrak{A}_d$-modules. It follows by the previous lemma that

$$C^\bullet_{\mathfrak{A}_{<d}}(S_d) \to C^\bullet_{\mathfrak{A}_{<d}}(C^\bullet_{\mathfrak{A}_d}(S)) \cong C^\bullet_{\mathfrak{A}_P}(S)$$

is a quasi-isomorphism.

Let $S_{d-1} \subset S_d$ be the subspace of functions constant along $A_{d-1}$ (and also along $A_d$). Repeating the argument above, we have $S_{d-1} \subset \Lambda^\bullet \mathfrak{A}_{d-1}^* \otimes S_d$ (inclusion in degree zero) induces an isomorphism $S_{d-1} \cong H^*_{\mathfrak{A}_{d-1}}(S_d)$ of $\mathfrak{A}_{<d-1}$-modules. This implies that

$$C^\bullet_{\mathfrak{A}_{<d-1}}(S_{d-1}) \to C^\bullet_{\mathfrak{A}_{<d-1}}(C^\bullet_{\mathfrak{A}_{d-1}}(S_d)) \cong C^\bullet_{\mathfrak{A}_{<d}}(S_d) \to C^\bullet_{\mathfrak{A}_P}(S)$$

is a quasi-isomorphism. Proceeding in this fashion shows that the inclusion (in degree zero) $S_1 \subset C^\bullet_{\mathfrak{A}_P}(S)$ is a quasi-isomorphism, where $S_1$ is the space of $A_P$-constant functions. But $S_1 = A'_\alpha \otimes S_{(2)}(\Gamma_{M_P} \backslash OK_P)$, which proves the lemma if $\lambda_P + \rho_P + \alpha \in {}^+\mathfrak{A}_P^*$.

Now suppose $r_i \leq 0$ for some $i$. Then the lemma asserts that $H^*_{\mathfrak{A}_P}(S) = 0$. From the Hochschild-Serre sequence for the pair $\mathfrak{A}_i \subset \mathfrak{A}_P$ one sees that it is enough to show that $H^*_{\mathfrak{A}_i}(S) = 0$. Now $H^*_{\mathfrak{A}_i}(S)$ is computed by the two-term complex

$$\Lambda^0 \mathfrak{A}_i^* \otimes S \to \Lambda^1 \mathfrak{A}_i^* \otimes S$$

with differential $f \mapsto \alpha_i \otimes R_i f$. If $r_i \leq 0$ it has zero kernel in $S$ and, with the same operator $I$ to invert the differential, it has zero cokernel. Then $H^*_{\mathfrak{A}_i}(S) = 0$. □

*Remark* 3.6. It is through this lemma, as used in the proof of Proposition 3.7 below, that the logarithmic modification makes its entry. As will become clear in that proof, it plays a crucial role in the main theorem.

3.6. *Stalk cohomology.* Let $\mathfrak{W}_P$ denote the set of weights of $A_P$ appearing in $H^*_{\mathfrak{N}_P}(E)$; it is completely described by Kostant's theorem ([7, III.3.1]). As a consequence (1) if $\mathbf{P} \subset \mathbf{Q}$ then the weights in $\mathfrak{W}_Q$ are restrictions to $A_Q \subset A_P$ of weights in $\mathfrak{W}_P$ and (2) $\mathfrak{W}_P$ is a subset of the weights appearing in $C^\bullet_{\mathfrak{N}_P}(E)$. Let $\mathfrak{W} = \mathfrak{W}_{P_0}$ for $\mathbf{P}_0$ minimal. Define two subsets of $\mathfrak{W}_P$ by

$$\mathfrak{W}_P(\overline{p}) = \{\, \alpha \in \mathfrak{W}_P \,|\, \alpha \in p_P + {}^+\mathfrak{A}_P^* \,\},$$
$$\mathfrak{V}_P(\lambda) = \{\, \alpha \in \mathfrak{W}_P \,|\, \alpha \in -\lambda_P - \rho_P + {}^+\mathfrak{A}_P^* \,\}.$$

By (1) and (2), $\mathfrak{W}(\overline{p}) = \mathfrak{V}(\lambda)$ if and only if $\mathfrak{W}_P(\overline{p}) = \mathfrak{V}_P(\lambda)$ for every $\mathbf{P}$.

PROPOSITION 3.7. *The inclusion*

$$\varinjlim_{t,O} W^\bullet_p(C_P(t,O), \mathbb{E}) \longrightarrow \varinjlim_{t,O} P\mathbf{S}^\bullet_{\lambda-\log}(C_P(t,O), \mathbb{E})^{N_P-\text{inv}}$$

*is an isomorphism in cohomology if* $-\lambda = \rho_0 + p$. *More generally, it is an isomorphism if* $\mathfrak{W}_P(\overline{p}) = \mathfrak{V}_P(\lambda)$.

*Proof.* From the discussion in 2.4 we know that there is an isomorphism

$$P\mathbf{S}^\bullet_{\lambda-\log}(C_P(t,O), \mathbb{E})^{N_P-\text{inv}}$$
$$\cong \bigoplus_\alpha C^\bullet_{(\mathfrak{M}_P, K_P)} \Big( C^\bullet_{\mathfrak{A}_P}(S_{\lambda+\alpha-\log}(\Gamma_{M_P} \backslash A_P(t) O K_P)) \otimes V^\bullet_\alpha \Big)$$

of complexes. Now Lemma 3.5 implies that the inclusion

$$A'_\alpha \otimes S_{(2)}(\Gamma_{M_P} \backslash O K_P) \otimes V^\bullet_\alpha \subset C^\bullet_{\mathfrak{A}_P}(S_{\lambda+\alpha-\log}(\Gamma_{M_P} \backslash A_P(t) O K_P))) \otimes V^\bullet_\alpha$$

is a quasi-isomorphism and induces an isomorphism of $(\mathfrak{M}_P, K_P)$-modules in cohomology. By Lemma 3.4 the inclusion

$$C^\bullet_{(\mathfrak{M}_P, K_P)}\big(A'_\alpha \otimes S_{(2)}(\Gamma_{M_P} \backslash O K_P) \otimes V^\bullet_\alpha\big) \longrightarrow$$
$$C^\bullet_{(\mathfrak{M}_P, K_P)}\big(S_{\lambda+\alpha-\log}(\Gamma_{M_P} \backslash A_P(t) O K_P)) \otimes V^\bullet_\alpha\big)$$

is a quasi-isomorphism for each $\alpha$. Recall from 2.5 that

$$W^\bullet_p(C_P(t,O)) \cong \bigoplus_\alpha C^\bullet_{(\mathfrak{M}_P, K_P)}\big(A_\alpha \otimes C^\infty_{\text{res}}(\Gamma_{M_P} \backslash O K_P) \otimes V^\bullet_\alpha\big).$$

Now
$$\varinjlim_O S_{(2)}(\Gamma_{M_P} \backslash O K_P) = \varinjlim_O C^\infty_{\text{res}}(\Gamma_{M_P} \backslash O K_P).$$

Indeed, both spaces are identified with the space $C_x$ of equivalence classes of pairs $(U, f)$ where $U$ is a neighbourhood of $\mathrm{pr}_{K_P}^{-1}(x)$ and $f$ is a smooth function on $U$, up to the equivalence $(U, f) \sim (U', f')$ if $f, f'$ agree on $U \cap U'$. (Note that $C^\bullet_{(\mathfrak{M}_P, K_P)}(C_x \otimes V_\alpha^\bullet)$ is identified with the complex $\Omega^\bullet(\mathbb{V}_\alpha^\bullet)_x$ of germs of smooth $\mathbb{V}_\alpha^\bullet$-valued forms at $x$.) Then

$$A_\alpha \otimes C_x \otimes V_\alpha^\bullet \subset A'_\alpha \otimes C_x \otimes V_\alpha^\bullet$$

induces an isomorphism of $(\mathfrak{M}_P, K_P)$-modules in cohomology precisely when $\mathfrak{V}_P(\lambda) = \mathfrak{W}_P(\overline{p})$. By Lemma 3.4 this implies that

$$\bigoplus_\alpha C^\bullet_{(\mathfrak{M}_P, K_P)}(A_\alpha \otimes C_x \otimes V_\alpha^\bullet) \subset \bigoplus_\alpha C^\bullet_{(\mathfrak{M}_P, K_P)}(A'_\alpha \otimes C_x \otimes V_\alpha^\bullet)$$

is a quasi-isomorphism under this condition. Up to the isomorphisms above, this is the map of the proposition. □

3.7. *The main theorem.* Let $\mathfrak{W}_P(\underline{p}) = \{\, \alpha \in \mathfrak{W}_P \,|\, \alpha \in p_P + \overline{{}^+\mathfrak{A}_P^*}\,\}$ and $\mathfrak{V}_P^+(\lambda) = \{\, \alpha \in \mathfrak{W}_P \,|\, \alpha \in -\lambda_P - \rho_P + \overline{{}^+\mathfrak{A}_P^*}\,\}$. (When **P** is maximal these differ from $\mathfrak{W}_P(\overline{p})$ and $\mathfrak{V}_P(\lambda)$ in at most one weight.) We are now in a position to state and prove the main theorem, which implies Theorem A and Corollary A in the introduction:

THEOREM 3.8. *The inclusions $\mathbf{W}^{\overline{p}}\mathbf{C}^\bullet(\mathbb{E}) \to \mathbf{S}^\bullet_{\lambda-\log}(\mathbb{E})$ and $\mathbf{W}^{\underline{p}}\mathbf{C}^\bullet(\mathbb{E}) \to \mathbf{S}^\bullet_{\lambda+\log}(\mathbb{E})$ are quasi-isomorphisms if $\lambda = -p - \rho_0$. More generally, $\mathbf{W}^{\overline{p}}\mathbf{C}^\bullet(\mathbb{E}) \to \mathbf{S}^\bullet_{\lambda-\log}(\mathbb{E})$ (resp. $\mathbf{W}^{\underline{p}}\mathbf{C}^\bullet(\mathbb{E}) \to \mathbf{S}^\bullet_{\lambda+\log}(\mathbb{E})$) is a quasi-isomorphism for any p and $\lambda$ for which $\mathfrak{W}(\overline{p}) = \mathfrak{V}(\lambda)$ (resp. $\mathfrak{W}(\underline{p}) = \mathfrak{V}^+(\lambda)$).*

*Proof.* I prove only the assertion about $\mathbf{S}^\bullet_{\lambda-\log}(\mathbb{E})$; the proof carries over to the other case (or one can use duality). By Propositions 3.2 and 3.7 the map $\mathbf{W}^{\overline{p}}\mathbf{C}^\bullet(\mathbb{E}) \to \mathbf{S}^\bullet_{\lambda-\log}(\mathbb{E})$ induces an isomorphism on stalk cohomology at any point $x$ on the boundary. Since this is clearly also true over $X$ the theorem follows. □

3.8. *Middle profiles.* For the upper and lower middle profiles, $\mu = -\overline{\rho_0}$ and $\nu = -\underline{\rho_0}$, $\mathfrak{W}(\mu) = \mathfrak{V}(0)$ and $\mathfrak{W}(\nu) = \mathfrak{V}^+(0)$. By Theorem 3.7,

COROLLARY 1. *The upper middle complex of sheaves $\mathbf{W}^\mu \mathbf{C}^\bullet(\mathbb{E})$ is quasi-isomorphic to $\mathbf{S}^\bullet_{(2)-\log}(\mathbb{E})$ and the lower middle complex of sheaves $\mathbf{W}^\nu \mathbf{C}^\bullet(\mathbb{E})$ is quasi-isomorphic to $\mathbf{S}^\bullet_{(2)+\log}(\mathbb{E})$.*

## 4. $(\mathfrak{G}, K)$-cohomology and $L^2$ cohomology

In this section Theorem B is proved using results of Franke's [11]. It would seem, *a priori*, that these are only valid for congruence $\Gamma$ because Franke uses Arthur's adelic truncation theory ([11, §5.3]). Truncation theory is available [23] for an arbitrary arithmetic group, so this assumption on $\Gamma$ is probably unnecessary. I have not, however, checked that the precise estimates needed in [11] are in [23], so the reader may wish to assume for this section (and for 5.1, 5.2) that $\Gamma$ is a congruence subgroup.

4.1. *Automorphic forms.* In computing $H^*_{\lambda-\log}(X, \mathbb{E})$, the module $S_{\lambda-\log}(\Gamma\backslash G)$ may sometimes be replaced by a submodule of automorphic forms (in the sense of [17]). Let $\mathcal{I} \subset \mathcal{Z}(\mathfrak{G})$ be the ideal annihilating $E^*$. By Franke's theorem ([11, Th. 15] or [28, 6.5]) (due in the $\mathbb{Q}$-rank one case to Casselman and Speh) the inclusion of the submodule $A_{\lambda-\log,\mathcal{I}}$ of $\mathcal{I}$-finite vectors in $S_{\lambda-\log}$ induces an isomorphism in $(\mathfrak{G}, K)$-cohomology when $\lambda$ is in the (open) dominant Weyl chamber $(\mathfrak{A}^*)^+$. Therefore for profiles $p$ such that $-p \in \rho_0 + (\mathfrak{A}^*)^+$ weighted cohomology can be computed using automorphic forms:

$$W^{\overline{p}}H^*(\Gamma, E) \simeq H^*_{(\mathfrak{G},K)}(A_{-p-\rho_0-\log,\mathcal{I}}(\Gamma\backslash G) \otimes E).$$

This includes the lower middle profile $\nu$. The profiles satisfying $p \in -\rho_0 + (\mathfrak{A}^*)^+$ (including upper middle) are dual to these.

4.2. $L^2$ *cohomology.* The inclusions

$$S_{(2)-\log}(\Gamma\backslash G) \subset S_{(2)}(\Gamma\backslash G) \subset S_{(2)+\log}(\Gamma\backslash G)$$

of $(\mathfrak{G}, K)$-modules induce maps on cohomology groups.

THEOREM 4.1. *If no proper parabolic subgroup of $G$ contains a Cartan subgroup of $K$ (in particular, if the ranks of $G$ and $K$ are equal) then the maps*

$$H^*_{(2)-\log}(X, \mathbb{E}) \to H^*_{(2)}(X, \mathbb{E}) \to H^*_{(2)+\log}(X, \mathbb{E})$$

*are isomorphisms.*

Note that by duality (1.9) it suffices to prove that the second map is an isomorphism. The proof (in §4.3) will use results due to Langlands [18], Franke [11] and Borel and Casselman [5].

4.3. *Proof of Theorem* 4.1. Let $S_{(2),d}$ and $S_{(2)+\log,d}$ be the discrete spectra of the spaces $S_{(2)}$ and $S_{(2)+\log}$. By results of Langlands' [18] and of Franke's ([11, Ths. 11, 12] or [28, 6.2, 6.3]) we know that in fact $S_{(2),d} = S_{(2)+\log,d}$. Consider the following diagram in which all maps are the obvious inclusions and each square commutes (notation as in 4.1):

$$\begin{array}{ccccc}
S_{(2),d} \otimes E & \xrightarrow{=} & S_{(2)+\log,d} \otimes E & \xleftarrow{a_1} & A_{(2)+\log,d,\mathcal{I}} \otimes E \\
i_1 \downarrow & & \downarrow i_2 & & \downarrow i_3 \\
S_{(2)} \otimes E & \xrightarrow{j} & S_{(2)+\log} \otimes E & \xleftarrow{a_2} & A_{(2)+\log,\mathcal{I}} \otimes E.
\end{array}$$

I want to show that $j$ induces an isomorphism in cohomology. By the result of Borel and Casselman ([5, Th. 4.5]), $i_1$ is an isomorphism in cohomology under the hypotheses of the theorem, so it will suffice to show that $i_2$ is also one. Franke has shown ([11, Th. 13]) that both $a_1$ and $a_2$ induce isomorphisms, so it will suffice to show that $i_3$ does so. For this it must be shown that $(A_{(2)+\log,\mathcal{I}}/A_{(2)+\log,d,\mathcal{I}}) \otimes E$ has no $(\mathfrak{G}, K)$-cohomology.

By [11, Th. 13] ([28, 4.4]), a constituent of $A_{(2)+\log,\mathcal{I}}/A_{(2)+\log,d,\mathcal{I}}$ is the summand fixed by a finite group (intertwining operators) in a (finite) direct sum of modules of the form

$$\Pi = \mathrm{Ind}_{(\mathfrak{P},K_P)}^{(\mathfrak{G},K)}(\pi \otimes D'_F).$$

The notation is as follows: $\mathbf{P}$ is a proper rational parabolic subgroup, $\pi$ is an irreducible subrepresentation in the space of cuspidal $L^2$ functions on $\Gamma_{M_P} \backslash M_P$ on which $\mathcal{Z}(\mathfrak{M}_P)$ acts by a character $\xi$, and $D'_F$ is the space of distributions on $i\mathfrak{A}_P^*$ supported on a certain finite set $F$ depending on $\xi$ (the exact definition of $F$ is unnecessary for the sequel). It suffices to show that $H^*_{(\mathfrak{G},K)}(\Pi \otimes E) = 0$. By Shapiro's lemma ([7, III.2.5]),

$$H^*_{(\mathfrak{G},K)}(\Pi \otimes E) \cong H^*_{(\mathfrak{P},K_P)}(\pi \otimes D'_F \otimes E).$$

The Hochschild-Serre spectral sequence for the pair $\mathfrak{N}_P \subset (\mathfrak{P}, K_P)$ has

$$E_2^{p,q} = H^p_{(\mathfrak{M}_P + \mathfrak{A}_P, K_P)}(\pi \otimes D'_F \otimes H^q_{\mathfrak{N}_P}(E)).$$

The spectral sequence degenerates at the $E_2$ term by a standard argument using Kostant's theorem ([7, III.3.4]). Let $W^{\mathbf{P}}$ be the set of minimal length representatives for $W/W_{\mathbf{P}}$, let $w \cdot \Lambda = w(\Lambda + \rho_P) - \rho_P$, and let $E^{\mathbf{M_P A_P}}_{w \cdot \Lambda}$ be the irreducible $\mathbf{M_P S_P}$-module with highest weight $w \cdot \Lambda$. Then the $E_2 = E_\infty$ term is the sum over $w \in W^{\mathbf{P}}$ of

$$H^{i-\ell(w)}_{(\mathfrak{M}_P + \mathfrak{A}_P, K_P)}(\pi \otimes D'_F \otimes E^{M_P A_P}_{w \cdot \Lambda}).$$

By a Künneth formula ([7, I.1.3]) this is equal to

$$\bigoplus_{p+q=i-\ell(w)} H^p_{(\mathfrak{M}_P, K_P)}(\pi \otimes E^{M_P}_{w \cdot \Lambda}) \otimes H^q_{\mathfrak{A}_R}(\mathbb{C}_{w \cdot \Lambda | \mathfrak{A}_P} \otimes D'_F).$$

LEMMA 4.2. $H^i_{\mathfrak{A}_P}(\mathbb{C}_\chi \otimes D'_\lambda) = 0$ if $\chi + \rho_P \neq -\lambda$ or $i > 0$ and $H^0_{\mathfrak{A}_P}(\mathbb{C}_{-\rho_P - \lambda} \otimes D'_\lambda) = \mathbb{C}$.

At this stage we can complete the proof of the theorem. Borel and Casselman have shown ([5, 1.8, 3.4]) that if $H^*_{(\mathfrak{M}_P, K_P)}(\pi \otimes E^{M_P}_{w \cdot \Lambda}) \neq 0$ then $P = \mathbf{P}(\mathbb{R})$ must contain a Cartan subgroup of $K$. (Briefly, the nonvanishing implies equality of the characters of $\mathcal{Z}(\mathfrak{M}_P)$ on $\pi$ and $E^{\mathbf{M_P}}_{w \cdot \Lambda}$, which implies that $E^{\mathbf{M_P}}_{w \cdot \Lambda}$ is isomorphic to its complex conjugate dual representation. This last fact is shown ([5, 1.8]) to imply the condition on $P$.) So, if no proper rational parabolic subgroup of $\mathbf{G}$ contains a Cartan subgroup of $K$, the contribution of nondiscrete Eisenstein series (i.e. of modules like $\Pi \otimes E$) must vanish. □

*Proof of Lemma* 4.2. Recall that the $\mathfrak{A}_P$-action on $D'$ (= distributions on $i\mathfrak{A}_P^*$) is given by

$$(\xi f)(i\eta) = \langle \xi, i\eta + \rho_P \rangle f(i\eta) \qquad (\xi \in \mathfrak{A}_P, f \in D', \eta \in \mathfrak{A}_P^*).$$

Let $\alpha_1, \ldots, \alpha_{d(P)}$ be the restrictions of the simple roots nontrivial on $A_P$. They give a basis of $\mathfrak{A}_P^*$ and so coordinates $x_1, \ldots, x_{d(P)}$ on $\mathfrak{A}_P^*$. The Lie algebra cohomology complex $C^\bullet(\mathfrak{A}_P, \mathbb{C}_\chi \otimes D'_\lambda)$ can therefore be identified ([26, 51.6]) with the topological tensor product of the $d(P)$ complexes $D'(i\mathbb{R}\alpha_j)_{\langle \check{\alpha}_j, \lambda \rangle} \longrightarrow D'(i\mathbb{R}\alpha_j)_{\langle \check{\alpha}_j, \lambda \rangle}$ with differential given by $f(ix) \mapsto (ix + \langle \check{\alpha}_j, \chi + \rho_P \rangle) f(ix)$. If $\langle \check{\alpha}_j, \chi + \rho_P \rangle = -\langle \check{\alpha}_j, \lambda \rangle$ this complex has cohomology $\mathbb{C}[0]$ and otherwise has zero cohomology. Since the spaces $D'(i\mathbb{R}\alpha_j)$ are nuclear, the lemma follows from standard properties of the topological tensor product [7, IX.6.1.3]). □

*Remark* 4.3. Something equivalent to the theorem is conjectured in [30] (first isomorphism in 8.4) for arbitrary $\mathbf{G}$, but this cannot be (by the above and [5, 4.6]).

## 5. Complements

5.1. *Trace of Hecke operators.* Suppose the ranks of $G$ and $K$ are equal. Then we have an isomorphism (Corollary A in the introduction)

$$W^\nu H^i(\Gamma, E) \simeq H^i_{(2)}(X, \mathbb{E}) \simeq W^\mu H^i(\Gamma, E).$$

When $X$ is Hermitian this follows from [20], [24] and [12, 23.2]. Starting from this equality, Goresky, Kottwitz, MacPherson ([15], [16], [14]) have given a topological proof of Arthur's formula [1] for the Lefschetz numbers of Hecke operators on the $L^2$ cohomology of $\mathbf{G}(\mathbb{Q}) \backslash \mathbf{G}(\mathbb{A}) / K_f K$. (Here $\mathbb{A}$= adeles, $\mathbb{A}_f$= finite adeles, and $K_f$ is a subgroup of $\mathbf{G}(\mathbb{A}_f)$.) Substituting the above isomorphism we see that their approach rederives Arthur's formula in the full generality of [1], i.e., under the equal-rank assumption. Even in the Hermitian case this is a simplification, since it avoids Looijenga's Hodge theory methods.

5.2. *The equal-rank case.* Suppose again that $X$ is equal-rank. Let $\hat{X}$ be a Satake compactification of $X$ in which all boundary components are equal-rank. Borel has conjectured that $H^i_{(2)}(X, \mathbb{E}) \simeq IH^i(\hat{X}, \mathbb{E})$. By Corollary B of the introduction this would imply $W^\mu H^i(\Gamma, E) \simeq W^\nu H^i(\Gamma, E) \simeq IH^i(\hat{X}, \mathbb{E})$. I expect that (for the quotient map $\Phi : \overline{X} \to \hat{X}$) both $R\Phi_* \mathbf{W}^\mu \mathbf{C}^\bullet(\mathbb{E})$ and $R\Phi_* \mathbf{W}^\nu \mathbf{C}(\mathbb{E})$ are quasi-isomorphic to the intersection complex on $\hat{X}$.

5.3. $L^2$ *Euler characteristic.* Suppose that $X$ is equal-rank or, more generally, that the hypothesis of Theorem 4.1 holds. Then Theorem 4.1, combined with the considerations of [12, 17.9], give a formula for the Euler characteristic of the $L^2$ cohomology of $\Gamma$ in terms of the usual Euler characteristic of Levi subgroups. (This generalizes the formula in the Hermitian case in [12] and [25].) For details and an application see [13, 6.12], [22].

5.4. *The map $H^*_\lambda \to H^*$ and a theorem of Borel's.* Let $n_\lambda$ be the smallest integer such that for $i \leq n_\lambda$ every weight $\alpha$ appearing in $H^i_\mathfrak{N}(E)$ satisfies $\alpha + \lambda + \rho_0 \in {}^+\mathfrak{A}^*$. Then $\mathbf{W}^{\overline{-\lambda-\rho_0}} \mathbf{C}^\bullet(\mathbb{E}) \longrightarrow \mathbf{W}^{-\infty} \mathbf{C}^\bullet(\mathbb{E}) \simeq Ri_* \mathbb{E}$ is a quasi-isomorphism in degrees $\bullet \leq i$. The condition also implies that $\mathbf{S}^\bullet_{\lambda-\log}(\mathbb{E}) \to \mathbf{S}^\bullet_\lambda(\mathbb{E})$ is a quasi-isomorphism for $i \leq n_\lambda$. Therefore the map $H^i_\lambda \to H^i$ is an isomorphism for $i \leq n_\lambda$. For $\lambda = 0$ this is [29, 3.20], which improved a result of [3]. In terms of the weighted cohomology sheaves the improvement is the difference between isomorphism and quasi-isomorphism of sheaves. I thank Professor Borel for suggesting this application.

The $G$-invariant forms on $G/K$ descend to $X$ and are well-known to be harmonic, so there is a map $I^*_G := \{\text{invariant forms}\} \to H^*(\Gamma, \mathbb{C})$. Since $I^*_G = H^*_{(\mathfrak{G},K)}(\mathbb{C})$ (here $\mathbb{C}$ is the trivial representation), there is an injection $I^*_G \hookrightarrow H^*_{(2)}(X, \mathbb{C})$, which is a surjection up to degree $m_G$ where $m_G + 1$ is the first degree in which a nontrivial unitary representation of $G$ may have $(\mathfrak{G}, K)$-cohomology (note that $m_G + 1$ is at least the real rank of $\mathbf{G}$). By the previous paragraph, it follows that $I^*_G \longrightarrow H^*(\Gamma, \mathbb{C})$ is an isomorphism for $* \leq \min(n_0, m_G)$. This is Borel's result [3].

5.5. Weighted cohomology can also be constructed with more sophisticated truncation procedures using either Kostant's theorem or the stratification (see remark in [12, 35.4]). The resulting groups would *not*, in general, be expressible as weighted $L^2$ cohomology.

5.6. *Rational structure.* The groups $W^{\overline{p}} H^i(\Gamma, E)$ have natural $\mathbb{Q}$-structures ([12, IV]) and (via the isomorphism $W^{\overline{p}} H^i(\Gamma, E) \cong H^i_{\lambda-\log}(X, \mathbb{E})$ of Theorem A) also have a decomposition by associate classes of parabolic subgroups ([11, §3]). Whether this decomposition is rational is not well understood. In the case of the full cohomology (i.e., very positive $\lambda$) this is proved by Franke ([11, Th. 20]) for $\mathbf{G}$ the restriction of scalars of $\mathrm{GL}(n)$ from a number field to $\mathbb{Q}$.

5.7. *$L^p$ cohomology.* The methods of Sections 2 and 3 also give isomorphisms between the (log-modified) $L^p$ cohomology of $X$ and weighted cohomology groups for the appropriate profile (see [21]).

5.8. *Generalizations.* It is not necessary to assume that $\Gamma$ is neat provided we work throughout with differential forms on manifolds with orbifold singularities. Furthermore, all the results hold more generally for reductive $\mathbf{G}$, and for local systems coming from representations $E$ which are multiplicity-free for the split centre $\mathbf{S_G}$ and such that the $\mathbf{S_G}$-isotypes are irreducible ([13, 2.8]). Finally, the main theorems extend to the adelic situation and, with the appropriate modifications, give Hecke equivariant isomorphisms (this is implicit in 5.1 above).

School of Mathematics, Institute for Advanced Study, Princeton, NJ

*Current address*: School of Mathematics, Tata Institute of Fundamental Research, Bombay, India

*E-mail address*: arvind@math.tifr.res.in

## References

[1] J. Arthur, The $L^2$-Lefschetz numbers of Hecke operators, Invent. Math. **97** (1989), 257–290.

[2] A. Borel, *Introduction aux Groupes Arithmétiques*, Hermann, Paris (1969).

[3] ______, Stable real cohomology of arithmetic groups, Ann. Sci. École Norm. Sup. **7** (1974), 235–272; II, *Manifolds and Lie Groups* (J. Hano et al., eds.), Prog. in Math. **14**, Birkhäuser (1981), 21–55.

[4] ______, Regularization theorems in Lie algebra cohomology, Applications, Duke Math. J. **50** (1983), 605–623.

[5] A. Borel and W. Casselman, $L^2$-cohomology of locally symmetric manifolds of finite volume, Duke Math. J. **50** (1983), 625–647.

[6] A. Borel and J-P. Serre, Corners and arithmetic groups, Comment. Math. Helv. **48** (1973), 436–491.

[7] A. Borel and N. Wallach, *Continuous Cohomology, Discrete Subroups and Representations of Reductive Groups*, Ann. of Math. Studies **94**, Princeton University Press, Princeton (1980).

[8] J.-L. Brylinski and J.-P. Labesse, Cohomologie d'intersection et fonctions $L$ de certaines variétés de Shimura, Ann. Sci. École Norm. Sup. **17** (1984), 361–412.

[9] W. A. Casselman, Introduction to the $L^2$-cohomology of arithmetic quotients of bounded symmetric domains, in *Complex Analytic Singularities*, Adv. Stud. in Pure Math. **8**, North-Holland, New York (1987), 69–93.

[10] ______, Introduction to the Schwartz space of $\Gamma\backslash G$, Can. J. Math. **41** (1989), 285–320.

[11] J. Franke, Harmonic analysis in weighted $L^2$-spaces, Ann. Sci. École Norm. Sup. (4) **31** (1998), 181–279.

[12] M. Goresky, G. Harder, and R. MacPherson, Weighted cohomology, Invent. Math. **116** (1994), 139-213.

[13] M. Goresky, G. Harder, R. MacPherson, and Arvind Nair, Local cohomology of Baily-Borel compactifications, preprint (1997).

[14] M. Goresky, R. Kottwitz, and R. MacPherson, Discrete series characters and the Lefschetz formula for Hecke operators, Duke Math. J. **89** (1997), 477–554; Correction to article: Duke Math. J. **92** (1998), 665–666.

[15] M. Goresky and R. MacPherson, Local contribution to the Lefschetz fixed point formula, Invent. Math. **111** (1993), 1–33.

[16] ———, Lefschetz numbers of Hecke correspondences, in *The Zeta Functions of Picard Modular Surfaces* (D. Ramakrishnan and R. P. Langlands, eds.), CRM Montreal (1992), 465–478.

[17] Harish-Chandra, *Automorphic Forms on Semisimple Lie Groups*, notes by J. G. M. Mars, Lecture Notes in Math. **62**, Springer-Verlag, Berlin (1968).

[18] R. P. Langlands, *On the Functional Equations Satisfied by Eisenstein Series*, Lecture Notes in Math. **544**, Springer-Verlag, Berlin (1976).

[19] ———, Modular forms and $l$-adic representations, in *Modular Functions of One Variable* II, Lecture Notes in Math. **349**, Springer-Verlag, Berlin (1973), 361–500.

[20] E. Looijenga, $L^2$ cohomology of locally symmetric varieties, Compositio Math. **67** (1988), 3–20.

[21] A. Nair, Weighted cohomology of arithmetic groups, Ph.D. thesis, University of Michigan (1996).

[22] ———, The $L^2$ Euler characteristic of arithmetic quotients, preprint (1998).

[23] M. S. Osborne and G. Warner, The Selberg trace formula II: Partition, reduction, truncation, Pac. J. Math. **106** (1983), 307–496.

[24] L. Saper and M. Stern, $L^2$-cohomology of arithmetic varieties, Ann. of Math. **132** (1990), 1–69.

[25] M. Stern, Lefschetz formulae for arithmetic varieties, Invent. Math. **115** (1994), 241–296.

[26] F. Trèves, *Topological Vector Spaces, Distributions and Kernels*, Academic Press, New York (1967).

[27] W. T. van Est, A generalization of the Cartan-Leray spectral sequence I, II, Indag. Math. **20** (1958), 399–413.

[28] J.-L. Waldspurger, Cohomologie des espaces de formes automorphes (d'après J. Franke), Séminaire Bourbaki 1995/96, Exposé No. 809, Astérisque **241** (1997), 139–156.

[29] S. Zucker, $L^2$ cohomology of warped products and arithmetic groups, Invent. Math. **70** (1982/83), 169–218.

[30] ———, $L^p$-cohomology and Satake compactifications, in *Prospects in Complex Geometry* (J. Noguchi and T. Ohsawa, eds.), Lecture Notes in Math. **1468**, 317–339, Springer-Verlag, Berlin (1991).